\documentclass[onefignum,onetabnum]{siamart250211}

\usepackage{amsmath, amssymb, amscd}
\usepackage{xcolor}
\usepackage{multirow}
\usepackage{subfigure} 
\usepackage{hyperref}

\usepackage{algorithm}
\usepackage{algpseudocode}

\usepackage{cprotect}  
 
\numberwithin{equation}{section}  

\usepackage{graphicx, epstopdf}

\headers{Spectral Method for Linearized Collision Operators}{T. Yin, Z. Cai, and Y. Wang}

\title{A fast Fourier spectral method for the linearized Boltzmann collision operator\thanks{Submitted to the editors DATE.
\funding{Tianai Yin was supported by the China Scholarship Council (CSC) (File No.202204890003). Zhenning Cai's work was supported by the Academic Research Fund of the Ministry of Education of Singapore under grant number A-8000965-00-00. This work of Yanli Wang is partially supported by the Foundation of the President of China Academy of Engineering Physics (YZJJZQ2022017), and the National Natural Science Foundation of China (Grant No. 12171026, U2230402, and 12031013).}}}

\author{Tianai Yin\thanks{School of Mathematical Science, Eastern Institute of Technology, Ningbo, Zhejiang 315200, China
  (\email{tyin@eitech.edu.cn}).}
\and Zhenning Cai\thanks{Department of Mathematics, National University of Singapore, Singapore 119076
  (\email{matcz@nus.edu.sg}).}
\and Yanli Wang\thanks{Beijing Computational Science Research Center, Beijing 100193, China
  (\email{ylwang@csrc.ac.cn}).}}

\begin{document} 
\maketitle

\newcommand{\image}{\textbf{\rm i}}

\newcommand{\bxi}{\boldsymbol{\xi}}
\newcommand{\bdeta}{\boldsymbol{\eta}}
\newcommand{\bdzeta}{\boldsymbol{\zeta}}
\newcommand{\bsigma}{\boldsymbol{\sigma}}
\newcommand{\bomega}{\boldsymbol{\omega}}
\newcommand{\bx}{\boldsymbol{x}}
\newcommand{\bv}{\boldsymbol{v}}
\newcommand{\bw}{\boldsymbol{w}}
\newcommand{\bu}{\boldsymbol{u}}
\newcommand{\bh}{\boldsymbol{h}}

\newcommand{\bz}{\boldsymbol{z}}
\newcommand{\by}{\boldsymbol{y}}
\newcommand{\bl}{\boldsymbol{\ell}}
\newcommand{\boldsymbolm}{\boldsymbol{m}}

\newcommand{\bg}{\boldsymbol{g}}
\newcommand{\bq}{\boldsymbol{q}}
\newcommand{\bj}{\boldsymbol{j}}
\newcommand{\bk}{\boldsymbol{k}}
\newcommand{\bp}{\boldsymbol{p}}
\newcommand{\hf}{\hat{f}}
\newcommand{\bbf}{\boldsymbol{f}}

\newcommand{\bZ}{\mathbb{Z}}
\newcommand{\bR}{\mathbb{R}}
\newcommand{\bS}{\mathbb{S}}

\newcommand{\mF}{\mathcal{F}}
\newcommand{\mD}{\mathcal{D}}
\newcommand{\mQ}{\mathcal{Q}}
\newcommand{\hQ}{\hat{Q}}
\newcommand{\hB}{\hat{B}}
\newcommand{\mB}{\mathcal{B}}

\newcommand{\mM}{\mathcal{M}}
\newcommand{\mL}{\mathcal{L}}
\newcommand{\dd}{{\rm d}}
\newcommand{\Kn}{\mathit{Kn}}
\newcommand{\Ma}{\mathit{Ma}}

\begin{abstract}
We introduce a fast Fourier spectral method to compute linearized collision operators of the Boltzmann equation for variable hard-sphere gases. While the state-of-the-art method provides a computational cost $O(M N^4 \log N)$, with $N$ being the number of modes in each direction and $M$ being the number of quadrature points on a hemisphere, our method reduces the cost to $O(N^4 \log N)$, removing the factor $M$, which could be large in our numerical tests. The method is applied in a numerical solver for the steady-state Boltzmann equation with quadratic collision operators. Numerical experiments for both spatially homogeneous and inhomogeneous Boltzmann equations have been carried out to test the accuracy and efficiency of our method.

\end{abstract}

\section{Introduction} \label{sec:introduction}
Rarefied gas dynamics is a branch of fluid dynamics that arises when the density of the gas is low (or when the characteristic length is small) and the rarefaction effects become more pronounced. It plays a vital role in astronautics and micro/nano flows \cite{shen2006rarefied}.
In astronautics, rarefied gas dynamics finds applications in many aerospace fields, including missiles, spacecraft, space shuttles, and space stations. It is also essential in micro-electro-mechanical systems, where the devices can be as tiny as a few microns. In these applications, rarefaction effects such as the gas-wall friction and the cold-to-hot heat transfer may dominate the gas flows. Simulations using rarefied gas models are therefore needed in such circumstances. Rarefied gas models are often built upon kinetic theories, which describe fluid states using distribution functions. In this work, we will focus on the Boltzmann equation, one of the fundamental kinetic models that incorporates the transport of gas molecules and their collisions into one integro-differential equation.

Solving the Boltzmann equation numerically has been a challenging task even with modern supercomputers, owing to its high dimensionality, complicated collision terms, and the gas-surface interaction that may ruin the regularity of the solution. In particular, the collision term that models the interaction between particles contains a five-dimensional integral, which results in high computational complexity if discretized using classical quadrature methods. At present, the mainstream numerical methods mainly fall into two categories: direct simulation of Monte Carlo (DSMC) and deterministic methods. The DSMC method \cite{MolecularBird} is a stochastic approach using simulation particles to model a cluster of gas molecules. In this method, the collision term is simulated according to the binary collision dynamics between the simulation particles, which avoids computation of the high-dimensional integral. It is a powerful and efficient solver for highly rarefied gas flows, but in the region where the gas is dense, its efficiency is weakened and the statistical noise may become strong.
Regarding the deterministic discretization of the collision term in the Boltzmann equation, the classical discrete velocity method \cite{DVM1964Broadwell,DVM1980,DVM1989} has evolved into spectral methods with remarkably higher efficiency, and in these methods,  the high-dimensional integrals are applied to the basis functions and can be precomputed. Spectral methods can also be categorized into two types: schemes on unbounded domains and schemes on truncated domains. The first category is a natural choice since the velocity domain is unbounded in the original Boltzmann equation, and different choices of basis functions lead to various methods, including the Hermite spectral method \cite{Kitzler2019,Hu2020}, the Burnett spectral method \cite{Gamba2018,Hu2020Burnett,CaiBurnett2020}, and the mapped Chebyshev method \cite{Hu2020Petrov,Hu2022}. However, due to the quadratic form of the collision term, the computational complexities of these methods are often close to the cube of the degrees of freedom in the velocity domain. Although some reduction of the complexity can be achieved by the sparsity of the coefficients, the computational time still grows quickly as the number of modes in the discretization increases. The second class of methods that apply the Fourier spectral method on truncated and periodized domains has a significantly lower computational cost. While the original version \cite{Pareschi1996} has the computational cost $O(N^6)$ ($N$ is the number of modes in each direction), more efficient schemes with computational complexity $O(N^{\alpha} \log N)$, $\alpha \leqslant 4$ have been developed for a general class of collision models \cite{mouhot2006fast,Wu2013, Hu2017}. Due to the truncation of the domain, some properties such as momentum and energy conservation are lost during the discretization, which needs to be fixed after every time step by a post-processing \cite{Gamba2009}.

Note that for the fast spectral methods with computational cost $O(N^{\alpha} \log N)$, the coefficient hidden in the big O often contains a factor $M$ being the number of quadrature points used to discretize an integral on the unit sphere $\mathbb{S}^2$. Although $\mathbb{S}^2$ is a two-dimensional manifold, the practical choice of $M$ can be much lower than $N^2$ \cite{mouhot2006fast, Hu2017}. In this work, we will also present numerical experiments to study the impact of this parameter on the accuracy of the numerical solutions.

Another approach to reducing the computational cost is to simplify the collision term in the Boltzmann equation. The most widely used method is the BGK-type approximations \cite{bhatnagar1954model,  shakhov1968generalization}, which replace the effect of binary collisions with an exponential relaxation to the equilibrium, with some fixes to correct some key physical quantities such as Prandtl numbers. These models are studied extensively in the literature, but experiments have also shown that they are not accurate enough in many cases \cite{Andries2002, Chen2015, Ambrus2020}. A better approximation is to linearize the collision operator about the local Maxwellian \cite{cai2018numerical}. Using the rotational invariance of the linearized collision operator, the numerical complexity can be reduced to $O(N^4)$ (without any large prefactors) if the Burnett spectral method is applied \cite{CaiBurnett2020}. However, for the Fourier spectral method, to the best of our knowledge, the computational cost of the linearized collision operator is no less than that of the binary collision operator \cite{Wu2014}.

The major contribution of this work is a novel numerical method that efficiently computes the linearized collision term discretized with the Fourier spectral method. With our approach, for variable hard-sphere gases, the computational complexity is $O(N^4 \log N)$ without the large prefactor $M$ that appeared in the fast spectral method for binary collision operators. A comparison between binary and linear collision operators applied to spatially homogeneous Boltzmann equation will be carried out to demonstrate the trade-off between the computational time and the numerical accuracy.

Applications of linearized collision operators to spatially inhomogeneous problems will also be represented in this paper. Instead of showing how linearized collision operators approximate the rarefied gas dynamics, we will apply our approach to a numerical solver of the steady-state nonlinear Boltzmann equation. Classically, the steady state of the Boltzmann equation can be obtained by long-time evolution of a certain initial state. Many efficient numerical methods to solve the time-dependent Boltzmann equation have been developed in the literature. Examples include asymptotic preserving schemes \cite{Bennoune2008, Dimarco2012}, the unified gas kinetic scheme \cite{XU2010}, the unified gas kinetic wave-particle method \cite{LIU2020108977}, the fast kinetic scheme \cite{Dimarco2018}, and dynamical low-rank methods \cite{lowrank2019}. In this work, we would like to apply iterative methods directly to find steady states. Existing methods of this type include the recently developed general synthetic iterative scheme \cite{Su2020, Zhu2021}. Our iterative scheme is based on the Newton method, which has quadratic convergence rate when the numerical solution during iteration is close to the exact solution. The original Newton method requires computing linearized Boltzmann collision operators. However, such linearized collision operators cannot be computed efficiently with our method since the linearization is performed about the solution at the previous step instead of the local Maxwellian. In our modified Newton method, we simply use the local Maxwellian to replace the solution at the previous step, so that our fast algorithm can be applied. Numerical experiments show that fast convergence can still be observed with this replacement.

In the rest part of the paper, we will introduce the Boltzmann equation and the existing Fourier spectral method \cite{Hu2017} for the collision operator in Section \ref{sec:preliminaries}. Our fast spectral method for the linearized collision operator is detailed in Section \ref{sec:numerical_scheme}. Section \ref{sec:num} presents numerical results showing the efficiency and accuracy of our method. Applications of our method to the steady-state Boltzmann equation are given in Section \ref{sec:application}. The paper ends with some concluding remarks in Section \ref{sec:conclusion}.

\section{Boltzmann collision operator and the Fourier spectral method} \label{sec:preliminaries}

\subsection{Boltzmann equation}
The Boltzmann equation is devised based on the kinetic description of gases by Ludwig Boltzmann in 1872, in which the fluid state is described by a probability distribution function (phase density function) $f(t,\bx,\bv)$, where $t$ is the time variable, and $\bx$ and $\bv$ represent the spatial location and velocity of gas molecules, respectively. For a single-species, monatomic gas without external forces, the Boltzmann equation has the following form
\begin{equation}
    \frac{\partial f(t,\bx,\bv)}{\partial t}
    +\bv\cdot\nabla_{\bx} f(t,\bx,\bv)
    =\mQ[f,f](t,\bx,\bv),
    \qquad t \in \mathbb{R}^+, \quad \bx \in \Omega, \quad \bv \in \mathbb{R}^3.
    \label{eq:Boltzmann_origin}
\end{equation}
The right-hand side $\mQ[f,f]$ describes binary collisions between gas molecules \cite{MolecularBird} and can be written as
\begin{equation}
    \mQ[f,f](t,\bx,\bv)
    =\int_{\bR^3} \int_{\bS^2}
    \mB(\bv-\bv_{\ast},\boldsymbol{\sigma})[f(t,\bx,\bv')f(t,\bx,\bv'_{\ast})-f(t,\bx,\bv)f(t,\bx,\bv_{\ast})]
    \,\dd \boldsymbol{\sigma} \,\dd \bv_{\ast},
    \label{eq:Q_origin}
\end{equation}
where the pre-collision velocities are $\bv$ and $\bv_{\ast}$, from which the post-collision velocities can be obtained by
\begin{equation}
    \bv' = \frac{\bv+\bv_{\ast}}{2}+\frac{|\bv-\bv_{\ast}|}{2}\boldsymbol{\sigma}, 
    \qquad \bv_{\ast}'=\frac{\bv+\bv_{\ast}}{2}-\frac{|\bv-\bv_{\ast}|}{2}\boldsymbol{\sigma}.
\label{eq:postcollision}
\end{equation}
The nonnegative collision kernel $\mB$ is determined by the intermolecular force between gas molecules. In this work, we are mainly concerned about the variable hard-sphere model, where $\mB$ is independent of $\boldsymbol{\sigma}$:
\begin{equation}
    \mB(\bv-\bv_{\ast},\boldsymbol{\sigma})= C |\bv-\bv_{\ast}|^{2(1-\omega)},
    \label{eq:mB}
\end{equation}
where the parameter $\omega \in [1/2, 1]$ is the viscosity index of the gas, and $C$ is a constant related to the diameter of the gas molecule. In particular, $\omega = 1$ and $\omega = 1/2$ correspond to the Maxwell molecules and the hard-sphere model, respectively. For simplicity, we will omit the parameter $\boldsymbol{\sigma}$ below and simplify write the collision kernel as $\mB(|\bv-\bv_{\ast}|)$.

As mentioned in the introduction, one approximation of the collision term $\mQ[f,f]$ is its linearized form about some local Maxwellian:
\begin{equation} \label{eq:linear_col}
\mathcal{L}[f] = \mQ[\mathcal{M},f] + \mQ[f,\mathcal{M}].
\end{equation}
Here the Maxwellian $\mathcal{M}$ has the expression
\begin{displaymath}
\mathcal{M}(\bv) = \frac{\rho}{(2\pi \theta)^{3/2}} \exp \left( -\frac{|\bv - \bu|^2}{2\theta} \right),
\end{displaymath}
where the quantities $\rho$, $\bu$ and $\theta$ refer to the density, velocity and temperature of the reference state, respectively. In many applications, the quantities $\rho$, $\bu$ and $\theta$ are chosen to be the moments of $f$:
\begin{displaymath}
\rho = \int_{\mathbb{R}^3} f \,\mathrm{d}\bv, \quad \bu = \frac{1}{\rho} \int_{\mathbb{R}^3} \bv f \,\mathrm{d}\bv, \quad \theta = \frac{1}{3\rho} \int_{\mathbb{R}^3} |\bv - \bu|^2 f \,\mathrm{d}\bv.
\end{displaymath}
In this work, we are also interested in some higher-order moments such as the heat flux $\bq$ and the pressure tensor $p_{ij}$, which are defined as
\begin{equation}
    \bq=\frac{1}{2}\int_{\bR^3}
    |\bv-\bu|^2(\bv-\bu)f\dd\bv,
    \qquad
    p_{ij}=\int_{\bR^3}(v_i-u_i)(v_j-u_j)f\dd \bv,\qquad i,j=1,2,3.
    \label{eq:qi_pij}
\end{equation}

In the literature, the numerical methods for $\mQ[f,f]$ has been well studied, but fewer methods can be found for $\mL[f]$. Before we introduce our method to compute $\mL[f]$ efficiently, we will first review a recent approach computing $\mQ[f,f]$ for arbitrary collision kernels.

\subsection{Fourier transform of the binary collision operator} \label{sec:Fourier}
Due to the complicated form and the high-dimensionality of the binary collision term \eqref{eq:Q_origin}, the most costly part of the Boltzmann numerical solver is the collision term. Currently, the Fourier spectral method is one of the most popular and efficient deterministic methods in the discretization of the collision term. The Fourier spectral method is generally developed based on the Fourier transform of the collision operator. In this section, we will review the Fourier transform of the binary collision operator and write the result in the form that matches the numerical method presented in \cite{Hu2017}, which can be applied to any elastic collision models.

Assume that $\hat{f}(\bxi)$ is the Fourier transform of the distribution function $f(\bv)$, satisfying
\begin{equation} \label{eq:fhat}
f(\bv) = \frac{1}{(2\pi)^3} \int_{\mathbb{R}^3} \hat{f}(\bxi) \mathrm{e}^{\image \bv \cdot \bxi} \,\mathrm{d}\bxi, \qquad \hat{f}(\bxi) = \int_{\mathbb{R}^3} f(\bv) \mathrm{e}^{-\image \bv \cdot \bxi} \,\mathrm{d} \bv.
\end{equation}
By definition, the Fourier transform of the collision term is
\begin{equation} \label{eq:Qhat}
\begin{aligned}
\widehat{\mQ}(\bxi) &= \int_{\mathbb{R}^3} \int_{\mathbb{R}^3} \int_{\mathbb{S}^2} \mathrm{e}^{-\image \bv \cdot \bxi} \mB(|\bv - \bv_{\ast}|) [f(\bv') f(\bv_{\ast}') - f(\bv) f(\bv_{\ast})] \,\dd \bsigma \,\dd \bv_{\ast} \,\dd \bv \\
&= \frac{1}{(2\pi)^6}\int_{\mathbb{R}^3} \int_{\mathbb{R}^3} \int_{\mathbb{R}^3} \int_{\mathbb{R}^3} \int_{\mathbb{S}^2}
\mB(|\bv - \bv_{\ast}|) \hat{f}(\bdeta) \hat{f}(\bdzeta) \\
& \hspace{50pt} \times (\mathrm{e}^{\image \bv' \cdot \bdeta}\mathrm{e}^{\image \bv_{\ast}' \cdot \bdzeta} - \mathrm{e}^{\image \bv \cdot \bdeta}\mathrm{e}^{\image \bv_{\ast} \cdot \bdzeta}) \mathrm{e}^{-\image \bv \cdot \bxi}  \,\dd \bsigma \,\dd \bv_{\ast} \,\dd \bv \,\dd \bdeta \,\dd \bdzeta.
\end{aligned}
\end{equation}
Further simplification requires introducing a set of new variables that is frequently used in manipulations of collision terms \cite{Kumar1966}:
\begin{equation} \label{eq:change_of_var}
\bg = \bv - \bv_{\ast}, \qquad \bg' = \bv' - \bv_{\ast}', \qquad \bh = \frac{\bv + \bv_{\ast}}{2}.
\end{equation}
It is not difficult to find
\begin{equation} \label{eq:gh}
\begin{gathered}
\bv = \bh + \frac{1}{2} \bg, \quad \bv_{\ast} = \bh - \frac{1}{2} \bg, \quad \bv' = \bh + \frac{1}{2} \bg', \quad \bv'_{\ast} = \bh - \frac{1}{2} \bg', \\
\bg' = |\bg| \bsigma, \quad \mathrm{d}\bg \,\mathrm{d}\bh = \mathrm{d}\bv \,\mathrm{d}\bv_{\ast}.
\end{gathered}
\end{equation}
These changes of variables enable us to further simply the expression of $\widehat{\mQ}$. The steps are:
\begin{enumerate}
\item Plug \eqref{eq:gh} into \eqref{eq:Qhat};
\item Integrate with respect to $\bh$;
\item Write $\bg$ in spherical coordinates $\bg = g \bomega$ and integrate with respect to $\bomega$.
\end{enumerate}
Here we omit the detailed calculations, and only show the final result:
\begin{equation} \label{eq:ft_binary}
\begin{aligned}
\widehat{\mQ}(\bxi) ={} & \frac{1}{\pi^2}  \int_0^{+\infty} \!\! \int_{\mathbb{S}_+^2} g^2 \mB(g) \operatorname{sinc} \left( \frac{g}{2} |\bxi|\right) I(g, \bxi,\bsigma) \dd\bsigma \,\dd g \\
&{} - \frac{2}{\pi} \int_{\mathbb{R}^3} \hat{f}(\bxi - \bdzeta) \hat{f}(\bdzeta)
  \left( \int_0^{+\infty} g^2 \mB(g) \operatorname{sinc} (g |\bdzeta|) \,\dd g\right) \,\dd\bdzeta,
\end{aligned}
\end{equation}
where
\begin{equation} \label{eq:I}
I(g, \bxi,\bsigma) = \int_{\mathbb{R}^3} \hat{f}(\bxi - \bdzeta) \hat{f}(\bdzeta) \left[ \cos \left( \frac{g \bsigma \cdot(\bxi - \bdzeta)}{2} \right) \cos \left( \frac{g \bsigma \cdot \bdzeta}{2} \right) + \sin \left( \frac{g \bsigma \cdot(\bxi - \bdzeta)}{2} \right) \sin \left( \frac{g \bsigma \cdot \bdzeta}{2} \right) \right] \dd \bdzeta.
\end{equation}
Note that this expression holds only formally, since in the second term of \eqref{eq:ft_binary}, the integral with respect to $g$ in the parentheses actually diverges. Nevertheless, the formula works when $\mB(\cdot)$ has a compact support, which is often used as an approximation of VHS models.

This expression shows the convolutional structure hidden in the collision term. In \eqref{eq:I}, when $g$ and $\bsigma$ are fixed, the integral with respect to $\bdzeta$ is the sum of two convolutions. In the second line, if the integral with respect to $g$ can be either analytically obtained or numerically precomputed, then this entire term is simply a convolution. Such a structure allows us efficiently compute the collision term using discrete Fourier transforms.

\subsection{Fourier spectral method for the binary collision operator}
\label{sec:fsm_binary}
To convert the continuous form \eqref{eq:ft_binary} to a numerical method. Several approximations must be carried out, including
\begin{itemize}
\item Since the Fourier spectral method can be applied only to periodic functions, the distribution functions are truncated and periodized so that $f \in L_{\mathrm{per}}^2([-L,L]^3)$.
\item The collision kernel $\mB(g)$ is truncated such that $\mB(g) = 0$ if $g > R$. For dealiasing purposes, $R$ should be selected such that $L \geqslant (3+\sqrt{2})R/4$.
\item The distribution function $f$ is discretized as
\begin{displaymath}
f(\bv) = \frac{1}{(2L)^3} \sum_{\bk = -N}^N \hat{f}_{\bk} \mathrm{e}^{\image \pi \bk \cdot \bv / L}, \qquad \hat{f}_{\bk} = \frac{1}{c_{\bk}} \left(\frac{L}{N} \right)^3 \sum_{\bl = -N}^{N-1} f\left( \frac{\bl L}{N} \right) \mathrm{e}^{-\image \pi \bk \cdot \bl / N},
\end{displaymath}
where $\bk = (k_1, k_2, k_3)$ and
\begin{displaymath}
\sum_{\bk = -N}^N = \sum_{k_1 = -N}^N \sum_{k_2 = -N}^N \sum_{k_3 = -N}^N, \quad c_{\bk} = c_{k_1} c_{k_2} c_{k_3}, \quad
c_k = \left\{ \begin{array}{@{}ll}
  2, & \text{if } k = \pm N, \\
  1, & \text{otherwise}.
\end{array} \right.
\end{displaymath}
Hereafter, we always assume that the multi-indices are cyclic, \text{i.e.} $(k_1 + 2n_1 N,k_2 + 2n_2 N,k_3 + 2n_3 N)$ are equivalent for all $(n_1, n_2, n_3) \in \mathbb{Z}^3$.
\item The integrals with respect to $g$ and $\bsigma$ are discretized using numerical quadratures:
\begin{equation} \label{eq:quadrature}
\int_0^R g^2 \psi(g) \,\mathrm{d}g \approx \sum_{j=1}^J w_j^g \psi(g_j), \qquad \int_{\mathbb{S}_+^2} \phi(\bsigma) \,\mathrm{d}\bsigma \approx 2\pi\sum_{m=1}^M w_m^{\sigma} \phi(\bsigma_m).
\end{equation}
Here $w_j^g$ and $w_m^{\sigma}$ are quadrature weights, and $g_j$ and $\bsigma_m$ are quadrature nodes.
\end{itemize}
With these approximations, the Fourier coefficients of the collision term are
\begin{equation} \label{eq:sm_binary}
\begin{aligned}
\widehat{\mQ}_{\bk} &= \frac{2\pi^2}{L^3} \sum_{j=1}^J \sum_{m=1}^M w_j^g w_m^{\sigma} \mathcal{B}(g_j) \operatorname{sinc} \left( \frac{\pi g_j |\bk|}{2L} \right) I_{j,\bk,m} \\
& \qquad - \frac{2\pi^2}{L^3} \sum_{\bl=-N}^{N-1} \hat{f}_{\bk-\bl} \hat{f}_{\bl} \left( \sum_{j=1}^J w_j^g \mB(g_j) \operatorname{sinc}\left( \frac{\pi g_j |\bl|}{L} \right) \right),
\end{aligned}
\end{equation}
where
\begin{displaymath}
\begin{aligned}
I_{j,\bk,m} = \sum_{\bl = -N}^{N-1} \hat{f}_{\bk-\bl} \hat{f}_{\bl} \Bigg[ &\cos \left( \frac{\pi g_j \bsigma_m \cdot(\bk - \bl)}{2L} \right) \cos \left( \frac{\pi g_j \bsigma_m \cdot \bl}{2L} \right) \\
& + \sin \left( \frac{\pi g_j \bsigma_m \cdot(\bk - \bl)}{2L} \right) \sin \left( \frac{\pi g_j \bsigma_m \cdot \bl}{2L} \right) \Bigg].
\end{aligned}
\end{displaymath}
In this formula, the summation over $\bl$ is a convolution, whose computational cost is $O(N^3 \log N)$. Taking into account the sums over $j$ and $m$, one can observe that the overall computational cost is $O(JMN^3 \log N)$. In the literature \cite{Hu2017, Jaiswal2019}, $J$ is usually chosen as $O(N)$, while $M$ is significantly smaller than $N^2$. As a result, the final computational cost is $O(M N^4 \log N)$.

In general, the Fourier spectral method does not preserve momentum and energy conservation, but is expected to preserve mass conservation, which means $\widehat{\mQ}_{\bk} = 0$ if $\bk = 0$. However, a straightforward application of \eqref{eq:sm_binary} does not preserve this property, which can be seen by observing that
\begin{equation} \label{eq:hatQ0}
\widehat{\mQ}_0 = \frac{2\pi^2}{L^3} \sum_{\bl=-N}^{N-1} \sum_{j=1}^J w_j^g \mB(g_j) \hat{f}_{-\bl} \hat{f}_{\bl} \left[ \sum_{m=1}^M w_m^{\sigma} \cos \left( \frac{\pi g_j \bsigma_m \cdot \bl}{L} \right) - \operatorname{sinc} \left( \frac{\pi g_j |\bl|}{L} \right) \right],
\end{equation}
is guaranteed to be zero only when the term in the square brackets equals zero. Two methods can be adopted to recover the mass conservation:
\begin{enumerate}
\item Simply set $\hat{Q}_0$ to zero.
\item Replace the sinc function in the second line of \eqref{eq:sm_binary} with the following approximation:
  \begin{equation} \label{eq:sinc}
    \operatorname{sinc} \left( \frac{\pi g_j |\bl|}{L} \right) = \frac{1}{2\pi} \int_{\mathbb{S}_+^2} \cos \left( \frac{\pi g_j \bsigma \cdot \bl}{L} \right) \,\dd \bsigma \approx \sum_{m=1}^M w_m^{\sigma} \cos \left( \frac{\pi g_j \bsigma_m \cdot \bl}{L} \right).
  \end{equation}
  This does not increase the computational cost as these values can be precomputed.
\end{enumerate}
Both methods will be tested in our numerical experiments.

\section{Fourier spectral method for the linearized collision operator} \label{sec:numerical_scheme}
In this section, we will study the numerical method for the linearized collision operator \eqref{eq:linear_col}. A naive implementation is to apply the algorithm for the binary collision operator to compute $Q[\mM,f]$ and $Q[f,\mM]$ separately, and add up the results. However, it is clear that this approach does not save any computational time and may be even more expensive than computing $Q[f,f]$ directly. Below, we will make use of the fact that one of the distribution functions is a Maxwellian, and introduce a scheme that can lower the time complexity to $O(N^4 \log N)$. Again, we begin with the the Fourier transform of the collision operator, and then reduce it to a numerical scheme by a series of approximations as described in the beginning of Section \ref{sec:fsm_binary}.

\subsection{Fourier transform of the linearized collision operator}
By definition, the linearized collision operator \eqref{eq:linear_col} can be written as $\mL(\bv) = \mL^+(\bv) - \mL^-(\bv)$,
where $\mL^+(\bv)$ and $\mL^-(\bv)$ are the gain term and the loss term, respectively, defined by
\begin{align*}
\mL^+(\bv) &= \int_{\mathbb{R}^3} \int_{\mathbb{S}^2} \mB(|\bv - \bv_{\ast}|) [\mM(\bv') f(\bv_{\ast}') + f(\bv') \mM(\bv_{\ast}')] \,\dd \bsigma \,\dd \bv_{\ast}, \\
\mL^-(\bv) &= \int_{\mathbb{R}^3} \int_{\mathbb{S}^2} \mB(|\bv - \bv_{\ast}|) [\mM(\bv) f(\bv_{\ast}) + f(\bv) \mM(\bv_{\ast})] \,\dd \bsigma \,\dd \bv_{\ast}.
\end{align*}
The Fourier transform of the loss term can be derived in the same way as the binary collision operator, and the result is
\begin{displaymath}
\widehat{\mL}^-(\bxi) = 
\frac{2}{\pi} \int_{\mathbb{R}^3} [\widehat{\mM}(\bxi - \bdzeta) \hat{f}(\bdzeta) + \hat{f}(\bxi - \bdzeta) \widehat{\mM}(\bdzeta)]
  \left( \int_0^{+\infty} g^2 \mB(g) \operatorname{sinc} (g |\bdzeta|) \,\dd g\right) \,\dd\bdzeta.
\end{displaymath}
Below we will focus on the gain term, which adopts a different approach.

Instead of expressing $\widehat{\mL}^+(\bxi)$ as an operator on $\hat{f}(\bxi)$, we now consider the function
\begin{equation} \label{eq:r}
r(\bv) = f(\bv) / \mM(\bv),
\end{equation}
and let $\hat{r}(\bxi)$ be the Fourier transform of $r(\bv)$.
Note that using $\hat{r}$ instead of $\hat{f}$ gives us a new formulation of $\widehat{\mL}^+(\bxi)$. The techniques used in our derivation again involve the change of variables \eqref{eq:change_of_var}. To avoid technical details, we omit the calculations and only present the final result:
\begin{equation} \label{eq:ft_linear}
\widehat{\mL}^+(\bxi) = \frac{4}{\pi} \int_0^{+\infty} \mB(g) \left[ \frac{\rho}{(4\pi\theta)^{3/2}} \exp \left( -\frac{g^2}{4\theta} \right) \right] \operatorname{sinc}\left( \frac{g|\bxi|}{2}\right) \int_{\mathbb{R}^3} H(g,\bv) \mathrm{e}^{-\image \bv \cdot \bxi} \,\dd\bv \,\dd g,
\end{equation}
where
\begin{equation} \label{eq:H}
H(g, \bv) = 8\pi^2 \mM^h(\bv) \int_{\mathbb{R}^3} \frac{f(\bv-\bw)}{\mM(\bv-\bw)} \delta\left( |\bw| - \frac{g}{2} \right) \dd \bw,
\end{equation}
and
\begin{displaymath}
\mM^h(\bh) = \frac{\rho}{(\pi\theta)^{3/2}} \exp\left( -\frac{|\bh-\bu|^2}{\theta} \right).
\end{displaymath}
Details of the derivation can be found in Supplementary Materials.
It can be observed that the integral over the hemisphere is no longer present in this result, so that the factor $M$ does not appear in the numerical method based on this expression.
The Fourier spectral method derived from \eqref{eq:ft_linear} will be introduced in the next subsection.

\subsection{Fourier spectral method for the linearized collision operator}
To convert \eqref{eq:ft_linear} to a numerical scheme, we again adopt the assumptions used in Section \ref{sec:fsm_binary}, including the truncation and periodization of the velocity domain, cut-off of the collision kernel $\mB(g)$, discretization of the distribution function $f$, and the approximation of the integral with respect to $g$ by numerical quadrature. We will reuse the notations in Section \ref{sec:fsm_binary} for consistency. Thus, the Fourier coefficients of the loss term has the same form as \eqref{eq:sm_binary}:
\begin{equation}
\label{eq:fsm_loss_term}
\widehat{\mL}_{\bk}^- = \frac{2\pi^2}{L^3} \sum_{\bl=-N}^N \left(\widehat{\mM}_{\bk-\bl} \hat{f}_{\bl} + \hat{f}_{\bk-\bl} \widehat{\mM}_{\bl} \right) \left( \sum_{j=1}^J w_j^g \mB(g_j) \operatorname{sinc}\left( \frac{\pi g_j |\bl|}{L} \right) \right),
\end{equation}
where
\begin{displaymath}
\widehat{\mM}_{\bl} = \widehat{\mM}\left( \frac{\pi \bl}{L} \right) = \rho \exp \left( -\frac{\image \pi}{L} \bl \cdot \bu - \frac{\pi^2 \theta}{L^2} |\bl|^2 \right), \qquad \bl = -N, \cdots, N.
\end{displaymath}
The computational cost of $\widehat{\mL}_{\bk}^-$ is again $O(N^3 \log N)$. Below we will focus on the computation of the gain term.

The general idea to approximate $\widehat{\mL}^+(\bxi)$ is to discretize the integral of $g$ using \eqref{eq:quadrature}. To compute the integrand we need to approximate $H(g_j, \bv)$ for each $g_j$ and apply the discrete Fourier transform. The definition of $H(g,\bv)$ \eqref{eq:H} shows that $H(g,\bv)$ is essentially $\mM^h(\bv)$ times a convolution. However, the convolution involves a Dirac delta function $\delta(|\bw| - g/2)$, which cannot be represented by interpolation on collocation points. To maintain the spectral accuracy, a natural idea to approximate the Dirac delta function is to interpolate its Fourier transform in the frequency space and then take the inverse discrete Fourier transform. However, the truncation of the frequency space may ruin the definition of $H(g,\bv)$. For easier calculation, we assume that the frequency domain is truncated to $B(0,\varrho)$, so that the Dirac delta function is approximated by
\begin{displaymath}
\begin{aligned}
\delta \left(|\bw| - \frac{g}{2} \right) \approx \Delta_{\varrho} (g,\bw) &= \frac{1}{(2\pi)^3} \int_{B(0,\varrho)} \pi g^2 \operatorname{sinc} \left( \frac{g|\bxi|}{2} \right) \mathrm{e}^{\image \bxi \cdot \bw} \,\dd \bxi \\
&= \frac{\varrho g}{2\pi |\bw|} \left[ \operatorname{sinc} \varrho \!\left( |\bw| - \frac{g}{2} \right) - \operatorname{sinc} \varrho \!\left( |\bw| + \frac{g}{2} \right) \right].
\end{aligned}
\end{displaymath}
This leads to the following approximation of $H(g,\bv)$:
\begin{displaymath}
H(g,\bv) \approx 8\pi^2 \mM^h(\bv) \int_{\mathbb{R}^3} \frac{f(\bv-\bw)}{\mM(\bv-\bw)} \Delta_{\varrho}(g, \bw) \,\dd\bw.
\end{displaymath}
Due to the slow decay of $\Delta_{\varrho}(g, \bw)$ with respect to $\bw$, when $f$ does not decay fast enough, the integral on the right-hand side of the above equation diverges.
For instance, when the distribution function $f$ has the form
\begin{displaymath}
f(\bv) = \frac{\rho_1}{(2\pi \theta_1)^{3/2}} \exp \left( -\frac{|\bv|^2}{2\theta_1} \right) + \frac{\rho_2}{(2\pi \theta_2)^{3/2}} \exp \left( -\frac{|\bv|^2}{2\theta_2} \right)
\end{displaymath}
with $\theta_1 < \theta_2$, the corresponding Maxwellian will have temperature between $\theta_1$ and $\theta_2$.
As a result, $f/\mathcal{M}$ increases to infinity as $|\bv| \rightarrow \infty$, and  thus the convolution between $f/\mathcal{M}$ and $\Delta_{\varrho}$ cannot be well defined.

Such an issue can be fixed by reformulating \eqref{eq:ft_linear} again. Note that
\begin{displaymath}
\begin{aligned}
& \frac{\rho}{(4\pi\theta)^{3/2}} \exp \left( -\frac{g^2}{4\theta} \right) H(g,\bv) \\
= {} & 8\pi^2 \mM^h(\bv) \int_{\mathbb{R}^3} \frac{\rho}{(4\pi\theta)^{3/2}} \exp \left( -\frac{g^2}{4\theta} \right) \frac{f(\bv-\bw)}{\mM(\bv-\bw)} \delta \left(|\bw| - \frac{g}{2} \right) \dd\bw \\
= {} & \pi^2 \mM^h(\bv) \int_{\mathbb{R}^3} \frac{\rho}{(\pi\theta)^{3/2}} \exp \left( -\frac{|\bw|^2}{\theta} \right) \frac{f(\bv-\bw)}{\mM(\bv-\bw)} \delta \left(|\bw| - \frac{g}{2} \right) \dd\bw.
\end{aligned}
\end{displaymath}
In this expression, even if $\delta(|\bw| - g/2)$ is replaced with $\Delta_{\varrho}(g,\bw)$, it is still well defined for all bounded $f$, due to the extra decaying term $\exp(-|\bw|^2/\theta)$. Thus, the final expression of $\widehat{\mL}^+(\bxi)$ to be converted to the numerical scheme is
\begin{equation} \label{eq:Lhat}
\widehat{\mL}^+(\bxi) = \frac{1}{2\pi} \int_0^{+\infty} g^2 \mB(g) \operatorname{sinc} \left( \frac{g|\bxi|}{2} \right) \left[ \int_{\mathbb{R}^3} \mM^h(\bv) \psi(g, \bv) \mathrm{e}^{-\image \bv \cdot \bxi} \,\dd\bv \right] \dd g,
\end{equation}
where
\begin{equation} \label{eq:psi}
\psi(g, \bv) = \int_{\mathbb{R}^3} \frac{f(\bv-\bw)}{\mM(\bv-\bw)} \frac{\rho}{(\pi\theta)^{3/2}} \exp \left( -\frac{|\bw|^2}{\theta} \right) \left[ \int_{\mathbb{R}^3} \operatorname{sinc} \left(\frac{g|\bdeta|}{2} \right) \mathrm{e}^{\image \bdeta \cdot \bw} \,\dd\bdeta\right] \dd\bw.
\end{equation}
We will now derive the numerical scheme based on this expression.

Let $\psi_{j,\bk}$ be the approximation of $\psi(g_j, \bk L/N)$. Then according to \eqref{eq:Lhat}, the Fourier coefficients of the gain term can be approximated by
\begin{equation} \label{eq:fsm_gain_term}
\widehat{\mL}_{\bk}^+ = \frac{L^3}{2\pi N^3} \sum_{j=1}^J w_j^g \mB(g_j)  \operatorname{sinc} \left( \frac{\pi g_j |\bk|}{2L} \right) \sum_{\bl = -N}^{N-1} \mM^h \left( \frac{\bl L}{N} \right) \psi_{j,\bl} \mathrm{e}^{-\image \pi \bk \cdot \bl / N},
\end{equation}
where the sum over $\bl$ can be computed by discrete Fourier transform, so that the total computational cost is $O(J N^3 \log N)$. The computation of $\psi_{j,\bk}$ follows \eqref{eq:psi}:
\begin{equation} \label{eq:psi_kl}
\begin{aligned}
\psi_{j,\bl} = \left( \frac{\pi}{N} \right)^3 \sum_{\bk = -N}^{N-1} r\left( \frac{(\bl - \bk) L}{N} \right)\frac{\rho}{(\pi\theta)^{3/2}} \exp \left( -\frac{|\bk|^2 L^2}{N^2 \theta} \right) \qquad \\
\times \left[ \sum_{\boldsymbol{m} = -N}^N \operatorname{sinc} \left( \frac{\pi g_j |\boldsymbol{m}|}{2L} \right) \mathrm{e}^{\image \pi \bk\cdot \boldsymbol{m} / N} \right].
\end{aligned}
\end{equation}
Here the sum over $\boldsymbol{m}$ can be computed via discrete Fourier transform, and the sum over $\bk$ is a convolution. Therefore, the total cost for computing all $\psi_{j,\bl}$ is also $O(J N^3 \log N)$. In practice, the term in the square brackets can be precomputed offline once $N$, $L$ and $g_j$ are determined, which may accelerate online computations by saving $J$ discrete Fourier transforms.

Compared with the binary collision term, whose computational cost is $O(M J N^3 \log N)$, using linearized collision operator divides the computational cost by a factor of $M$, which will lead to a remarkable reduction of the computational time. However, a direct application of this approach may cause significant loss of accuracy. We will start a new subsection to explain the reason.

\subsection{The accuracy issue}
The accuracy issue lies in the combination of \eqref{eq:r} and the discrete Fourier transform. Here we demonstrate the issue using a toy example. Consider the following sequence of operations:
\begin{displaymath}
\scalebox{.94}{%
$\begin{CD}
\boxed{f(\bv) := \mM(\bv)} @>\texttt{FFT}>> \boxed{\hspace{29pt}\hat{f}(\bv)\hspace{29pt}} @>\texttt{iFFT}>> \boxed{g(\bv)} @>\texttt{Div}>> \boxed{r(\bv) := g(\bv) / \mM(\bv)} \\
@. @. @. @V\texttt{FFT}VV \\
@. \boxed{q(\bv) := p(\bv) \mM(\bv)} @<\texttt{Mul}<< \boxed{p(\bv)} @<\texttt{iFFT}<< \boxed{\hspace{32pt}\hat{r}(\bv)\hspace{32pt}}
\end{CD}$}
\end{displaymath}
Under exact arithmetic, we should have
\begin{displaymath}
f(\bv) = g(\bv) = q(\bv) = \mM(\bv), \qquad r(\bv) = p(\bv) = 1.
\end{displaymath}
However, when these operations are implemented with floating-point numbers, some of these equalities are severely violated. As an example, we provide below the MATLAB code that carries out these operations:
\begin{center}
\smallskip
\cprotect\framebox{%
  \begin{minipage}{.9\linewidth}
    \begin{verbatim}
N = 16; L = 7.5;
x = (-N:N-1) * (L/N); [X, Y, Z] = meshgrid(x, x, x);
M = fftshift(1/(2*pi)^(3/2) * exp(-(X.^2+Y.^2+Z.^2)/2));
f = M; g = ifftn(fftn(f)); r = g./M;
p = ifftn(fftn(r)); q = p .* M;
fprintf("Max norm of f-g: %e\n", max(abs(f-g), [], 'all'));
fprintf("Max norm of f-q: %e\n", max(abs(f-q), [], 'all'));
fprintf("Max norm of r-1: %e\n", max(abs(r-1), [], 'all'));\end{verbatim}
  \end{minipage}
}
\smallskip
\end{center}
In this example, the parameters are
\begin{displaymath}
N = 16, \quad L = 7.5, \quad \rho = \theta = 1, \quad \bu = 0,
\end{displaymath}
and the result of our test is
\begin{center}
\smallskip
\cprotect\framebox{%
  \begin{minipage}{.9\linewidth}
    \begin{verbatim}
Max norm of f-g: 1.387779e-17
Max norm of f-q: 8.190679e+00
Max norm of r-1: 6.299552e+19\end{verbatim}
  \end{minipage}
}
\smallskip
\end{center}
It can be seen that although $f = g$ holds approximately, the values of $r$ and $q$ are way off. The reason is purely the round-off error of double-precision floating-point numbers.

The loss of accuracy due to round-off errors in algorithms using the fast Fourier transform has been noticed in several previous works \cite{Keich2005, Wilson2016, Wilson2017}. In general, although the difference between \texttt{f} and \texttt{g} looks small, the pointwise relative error could be large. For instance, in the example above, the values of \texttt{f(16,16,16)} and \texttt{g(16,16,16)} are \texttt{3.9486e-34} and \texttt{3.0542e-19}, respectively. Thus, the value of \texttt{r(16,16,16)}, which is \texttt{g(16,16,16)} divided by \texttt{f(16,16,16)}, becomes \texttt{7.7349e+14}, far away from the unity.

Our fast Fourier spectral method does not perform exactly the same transforms, but some operations resemble this simple example.
For instance, the computation of the convolutions in the loss term requires both forward and backward Fourier transforms of the distribution function, similar to the operations $f \rightarrow \hat{f} \rightarrow g$, and the gain term divides $f$ by $\mathcal{M}$, uses Fourier transforms to compute the convolution to get $\psi(g,\bv)$, and then multiplies $\psi(g,\bv)$ by $\mathcal{M}^h(\bv)$, which is analogous to the operations $g \rightarrow r \rightarrow \hat{r} \rightarrow p \rightarrow q$. Since the actual computation mixes up the gain term and the loss term,
the example above shows that our approximation of $r(\cdot)$ may contain large errors, especially in the tail of the distribution function. But even so, the computation following \eqref{eq:fsm_gain_term} and \eqref{eq:psi_kl} \emph{strictly} can still achieve high accuracy due to the Maxwellian $\mM^h(\bv)$ in \eqref{eq:Lhat} (or $\mM^h(\bl L/N)$ in \eqref{eq:fsm_gain_term}) that suppresses the function values for large $|\bv|$. Here by ``\emph{strictly}'', we mean that the convolution, \textit{i.e.}, the sum over $\bk$ in \eqref{eq:psi_kl}, is taken using direct summation without fast Fourier transforms. This method is referred to as ``naive convolution'' below. In fact, the step that really destroys the accuracy is the use of fast Fourier transform to compute this convolution, which requires a discrete Fourier transform on $r$. Such an issue has been studied in \cite{Wilson2016}, which shows that the naive convolution of two positive vectors can be far more accurate than convolution using fast Fourier transforms in terms of maximum relative errors. The reason is that the Fourier transforms may lead to catastrophic cancellations, which does not exist in the naive convolution. In our case, both functions in the convolution are close to positive functions, and the above statements still apply.

\subsection{Practical implementation} \label{sec:implementation}
In this section, we propose three solutions to tackle the difficulty raised in the previous section, including
\begin{itemize}
\item applying \emph{accurate FFT convolution (aFFT-C)} introduced in \cite{Wilson2016};
\item using floating-point numbers of higher precisions;
\item truncating the distribution function $f(\bv)$.
\end{itemize}
These methods will be discussed separately in the following paragraphs.

In the first method, we note that the general idea of the aFFT-C algorithm \cite{Wilson2016} is to decompose both functions in the convolution and convolve all the pairs. Each convolution is again computed with FFTs. The relative error is controlled by adopting the decompositions appropriately such that the support of each convolution does not spread too widely, allowing manual elimination of the numerical error outside its support. According to our test, this method indeed reduces the numerical error, but the computational cost is multiplied by the number of pairs, which may cancel out the factor $M$ we have saved compared with the computation of binary collision operators. This method is therefore not adopted in our experiments.

The second possible approach is to use floating-point numbers that are more accurate than double precision. In our experiments, we find that the quadruple-precision floating-point (\texttt{binary128}) format, which has approximately 34 significant decimal figures, is sufficient for practical use. The issue of this approach is that very few CPUs have native support of \texttt{binary128} numbers, and the speed of discrete Fourier transforms based on software implementations of \texttt{binary128} by math libraries such as GCC's libquadmath is usually significantly lower than those for double-precision floating-point numbers. We believe this method is practical on IBM's POWER9 processors, which have native hardware support for \texttt{binary128}.

The third idea is to replace the function $f$ by
\begin{equation} \label{eq:cutoff}
\tilde{f}(\bv) = \left\{ \begin{array}{@{}ll}
f(\bv), & \text{if } \rho^{-1} \mM(\bv) > \varepsilon, \\
0, & \text{if } \rho^{-1} \mM(\bv) < \varepsilon.
\end{array} \right.
\end{equation}
The purpose is to avoid excessively large values in $r(\bv)$, so that the precision of the convolution is controllable. The rationale behind this approach is the boundedness of the gain term in the collision operator \cite{Alonso2010}:
\begin{align*}
\|\mQ^+[f,g]\|_{L^s(\mathbb{R}^3)} \lesssim \left( \int_{\mathbb{R}^3} (1 + |\bv|^{2p\omega}) |f(\bv)|^p \,\dd \bv \right)^{1/p}\left( \int_{\mathbb{R}^3} (1 + |\bv|^{2q\omega}) |g(\bv)|^q \,\dd \bv \right)^{1/q}, \\ \forall p,q,s \in [1,+\infty] \text{ satisfying } 1/p + 1/q = 1 + 1/s.
\end{align*}
Such a result indicates that a small perturbation of $f$ only results in a small perturbation in the gain term $\mL^+[f]$. If $\mM(\bv)$ is the Maxwellian with the same density, velocity and temperature as $f(\bv)$, in most cases, for a $\bv$ such that $\mM(\bv)$ is small, $f(\bv)$ also has a small value, meaning that the difference between $f$ and $\tilde{f}$ is expected to be small when $\varepsilon$ is small. In our experiments, we choose $\varepsilon$ to be $10^{-9}$ for double-precision floating-point numbers. This method does not increase the computational cost and can still preserve good numerical accuracy. The only disadvantage is that the numerical error cannot be reduced to machine epsilon due to the cutoff \eqref{eq:cutoff}.

\subsection{Summary of the algorithm}
\label{sec:summary_of_algorithm}
As a summary, we rewrite our algorithm as a pseudocode in Algorithm \ref{algo:collision}, where both the input and output of the function \textsc{LinearizedCollisionOperator} are point values on collocation points, and we have rescaled the coefficients appropriately to match the definitions of the \textsc{FFT} and \textsc{InvFFT} in lines 29 to 34. The parameter $\textit{cutoff}$ of the function \textsc{LinearizedCollisionOperator} indicates whether the truncation \eqref{eq:cutoff} is needed. According to the discussion in Section \ref{sec:implementation}, when the datatype \texttt{binary128} is used, the parameter $\mathit{cutoff}$ in line 2 can be set to $\mathit{false}$, while $\mathit{cutoff}$ needs to be $\mathit{true}$ when double-precision numbers are used. Due to the slowness of \texttt{binary128} numbers, in our tests, arrays will be converted to the \texttt{binary128} type only when computing the convolutions, \textit{i.e.}, lines 13, 17 and 18 in Algorithm \ref{algo:collision}.

\begin{algorithm}[!ht]
\begin{algorithmic}[1]
\State \textbf{Precomputation:} For all $j = 1,\cdots,J$ and $\bk = -N, \cdots, N$, compute
\begin{gather*}
s_{j,\bk} \gets \operatorname{sinc} \left( \frac{\pi g_j |\bk|}{2L} \right), \qquad
\phi_{j,\bk} \gets 32\pi^2 w_j^g \mB(g_j) s_{j,\bk}, \\ \varphi_{j,\bk} \gets \textsc{InvFFT}(s_{j,\bk}) , \qquad \omega_{\bk} \gets 16\pi^2 \sum_{j=1}^J w_j^g \mB(g_j) \operatorname{sinc}\left( \frac{\pi g_j |\bk|}{L} \right)
\end{gather*}

\hrule height .2pt
\medskip

\Function{CollisionOperator}{$\rho, \bu, \theta, f_{\bk}, \mathit{cutoff}$}
\For{$\bk = -N, \cdots, N$}
\State $\bv_{\bk} \gets \bk L / N$, $M_{\bk} \gets \dfrac{\rho}{(2\pi\theta)^{3/2}} \exp \left( -\dfrac{|\bv_{\bk} - \bu|^2}{2\theta} \right)$
\State $M^g_{\bk} \gets \dfrac{\rho}{(\pi\theta)^{3/2}} \exp \left( -\dfrac{|\bv_{\bk}|^2}{\theta} \right)$,  $M^h_{\bk} = \dfrac{\rho}{(\pi\theta)^{3/2}} \exp \left( -\dfrac{|\bv_{\bk} - \bu|^2}{\theta} \right)$
\If{$\textit{cutoff}$ and $M_{\bk} / \rho < 10^{-9}$}
\State $r_{\bk} \gets 0$
\Else
\State $r_{\bk} \gets f_{\bk} / M_{\bk}$
\EndIf
\EndFor
\medskip
\State {\color{gray} // Compute the gain term}
\State $\hat{r}_{\bk} \gets \textsc{FFT}(r_{\bk})$
\State $\widehat{\mL}_{\bk}^+ \gets 0$,  $\forall \bk = -N, \cdots, N$
\For{$j = 1,\cdots, J$}
\State $u_{\bk} \gets M^g_{\bk} \varphi_{j,\bk}$, $\forall \bk = -N, \cdots, N-1$
\State $\hat{u}_{\bk} \gets \textsc{FFT}(u_{\bk})$
\State $\psi_{\bk} \gets \textsc{InvFFT}(\hat{u}_{\bk} \hat{r}_{\bk})$
\State $\psi_{\bk} \gets M_{\bk}^h \psi_{\bk}$ for all $\bk = -N,\cdots,N-1$
\State $\hat{\psi}_{\bk} = \textsc{FFT}(\psi_{\bk})$
\State $\widehat{\mL}_{\bk}^+ \gets \widehat{\mL}_{\bk}^+ + \phi_{j,\bk} \hat{\psi}_{\bk}$, $\forall \bk = -N, \cdots, N$
\EndFor
\State $\mL_{\bk}^+ \gets \textsc{InvFFT}(\widehat{\mL}_{\bk})$
\medskip
\State {\color{gray} // Add the loss term}
\State $\hat{f}_{\bk} \gets \textsc{FFT}(f_{\bk})$,  $\hat{M}_{\bk} \gets \textsc{FFT}(M_{\bk})$
\State $u_{\bk} \gets \textsc{InvFFT}(\hat{f}_{\bk} \omega_{\bk})$, $w_{\bk} \gets\textsc{InvFFT}( \hat{M}_{\bk} \omega_{\bk})$ \State $\mL_{\bk} \gets \mL_{\bk}^+ - u_{\bk} M_{\bk} - w_{\bk} f_{\bk}$, $\forall \bk = -N,\cdots,N$
\EndFunction
\medskip
\hrule height .2pt
\medskip
\Function{FFT}{$f_{\bk}$}
\State \Return $\displaystyle \hat{f}_{\bk} = \frac{1}{c_{\bk}} \sum_{\bl=-N}^{N-1} f_{\bl} \mathrm{e}^{-\image \pi \bk \cdot \bl / N}$, $\forall \bk = -N, \cdots, N-1$
\EndFunction
\medskip
\hrule height .2pt
\medskip
\Function{InvFFT}{$f_{\bk}$}
\State \Return $\displaystyle f_{\bk} = \frac{1}{(2N)^3} \sum_{\bl=-N}^N \hat{f}_{\bl} \mathrm{e}^{\image \pi \bk \cdot \bl / N}$, $\forall \bk = -N, \cdots, N-1$
\EndFunction
\end{algorithmic}
\caption{Linearized collision term}
\label{algo:collision}
\end{algorithm}

\clearpage

\section{Numerical results}
\label{sec:num}
\begin{figure}[!ht]
\centering
\subfigure[$f^{(1)}(v_1, v_2, 0)$]{
\includegraphics[scale=.5]{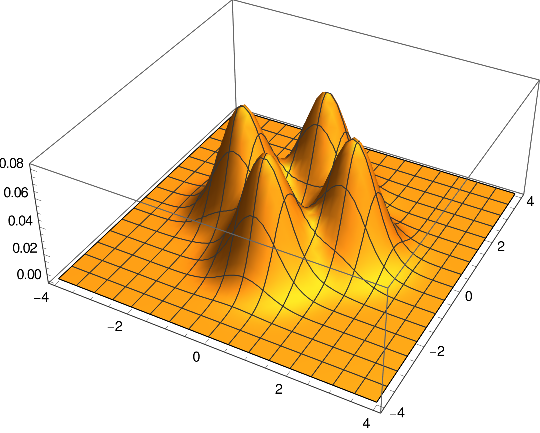}
} \quad
\subfigure[$f^{(2)}(v_1, v_2, 0)$]{
\includegraphics[scale=.5]{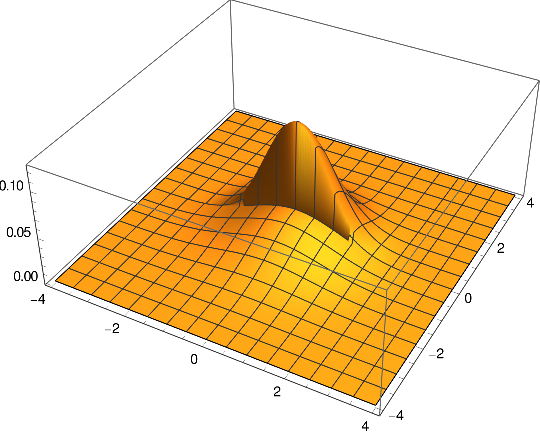}
}
\caption{Distribution functions $f^{(1)}$ and $f^{(2)}$ at $v_3 = 0$}
\label{fig:f1f2}
\end{figure}

In this section,  we will carry out some numerical tests to verify the validity of the Fourier spectral method for the linearized collision operator. In our experiments, we will test out two collision kernels:
\begin{equation} \label{eq:col_ker}
\mB_1(g) \equiv \frac{1}{4\pi}, \qquad \mB_2(g) = \frac{1}{4\pi} g^{0.56},
\end{equation}
where $\mB_1$ corresponds to Maxwell molecules, and $\mB_2$ corresponds to a gas with viscosity index $0.72$, which is often adopted for the argon gas. Two functions will be used frequently in our tests. The first one is a smooth distribution function:
\begin{displaymath}
\begin{aligned}
f^{(1)}(\bv) = \frac{1}{4\pi^{3/2}} & \left[ \exp \left( -\frac{(v_1+u)^2 + v_2^2 + v_3^2}{2\vartheta} \right) + \exp \left( -\frac{(v_1-u)^2 + v_2^2 + v_3^2}{2\vartheta} \right) \right. \\
& + \left. \exp \left( -\frac{v_1^2 + (v_2+u)^2 + v_3^2}{2\vartheta} \right) + \exp \left( -\frac{v_1^2 + (v_2-u)^2 + v_3^2}{2\vartheta} \right) \right],
\end{aligned}
\end{displaymath}
where $u = \sqrt{2}$ and $\vartheta = 1/3$. The second function is discontinuous:
\begin{displaymath}
f^{(2)}(\bv) = \left\{ \begin{array}{@{}ll}
\dfrac{\sqrt[4]{2}(2-\sqrt{2})}{\pi^{3/2}} \exp \left( -\dfrac{|\bv|^2}{\sqrt{2}} \right), & \text{if } v_1 > 0, \\[10pt]
\dfrac{\sqrt[4]{2}(2-\sqrt{2})}{4\pi^{3/2}} \exp \left( -\dfrac{|\bv|^2}{2\sqrt{2}} \right), & \text{if } v_1 < 0.
\end{array} \right.
\end{displaymath}
Both functions are designed such that the corresponding velocity is zero and the temperature is $1$. Figure \ref{fig:f1f2} shows the cross-section of these two functions at $v_3 = 0$.

\subsection{Accuracy test}
We first carry out numerical experiments to confirm the spectral accuracy of our method. Since the exact expression of $\mL[f]$ is generally unknown, we adopt two different methods to compute the collision term, and compare the results. The first method is Algorithm \ref{algo:collision} introduced in Section \ref{sec:summary_of_algorithm}, and the second is to apply the algorithm for quadratic collision operators introduced in \Cref{sec:fsm_binary} to compute $\mQ[f,\mM] + \mQ[\mM,f]$. In both methods, the integral over $g$ is approximated by the Jacobi-Gauss-Radau quadrature, and the number of quadrature point $J$ is chosen such that the numerical quadrature is exact for all polynomials of degrees up to $2N$. Numerical quadrature on the sphere, required by the quadratic collision operator, is also chosen such that all spherical harmonics of degree less than or equal to $2N$ can be exactly integrated. This can be achieved by using Lebedev quadrature nodes and weights, which offers formulas of orders up to 131. No conservation fix is imposed for both methods.

In our experiments, two test cases are considered:
\begin{description}
  \item[Case 1] $\mB = \mB_1$ and $f = f^{(1)}$;
  \item[Case 2] $\mB = \mB_2$ and $f = f^{(2)}$.
\end{description}
Several different values of $R$ (cutoff parameter of $g$) and $N$ (cutoff parameter for Fourier coefficients) are studied, and we present the $L^2$ differences in Tables \ref{tab:L2diff1} and \ref{tab:L2diff2}. In general, these two methods compute the loss term in exactly the same way, and therefore the difference lies only in the gain term. Three parameters may affect the difference between two methods:
\begin{enumerate}
\item truncation of the velocity domain: larger $L$ (hence larger $R$) leads to smaller difference, since our method relies on the equality \eqref{eq:ft_linear}, which holds exactly only when $L$ is infinity, while the method for the binary collision operator can be applied to periodized collision operators \cite{Hu2017}.
\item discretization of the velocity domain: larger $N$ leads to smaller difference, since our method is based on the Fourier coefficients of $r(\cdot)$ instead of $f(\cdot)$, meaning that the interpolations of the distribution functions on the physical space are essentially different, but they will converge to the same limit as $N \rightarrow \infty$.
\item The parameter for the cutoff of $f$, which is fixed to be $10^{-9}$ in our experiments using double-precision floating-point numbers.
\end{enumerate}
All the numbers in Table \ref{tab:L2diff1} and Table \ref{tab:L2diff2} are not affected by the cutoff of $f$, and this has been confirmed by experiments using quadruple-precision numbers without cutoff. For fixed $R$, the difference between two results decreases as $N$ increases, until the error caused by $R$ dominates. For sufficiently large $R$, the fast reduction of the error indicates the spectral accuracy of the method. If we fix $N$ and increase $R$, the error may decrease first if $N$ is sufficiently large so that the error caused by the truncation of the velocity domain is dominant when $R$ is small. Once $R$ is sufficiently large, decreasing $N$ will lead to a coarser grid, and thus gives greater differences between the two methods.
\begin{table}[!ht]
\centering
\caption{$L^2$ difference between two methods computing the linearized collision term (Case 1)}
\label{tab:L2diff1}
\scalebox{.9}{%
\begin{tabular}{c|cccc}
& $N = 4$ & $N = 8$ & $N = 16$ & $N = 32$ \\
\hline
$R = 4$ & $1.71 \times 10^{-2}$ & $1.35 \times 10^{-4}$ & $1.17 \times 10^{-4}$ & $1.17 \times 10^{-4}$ \\
$R = 6$ & $9.82 \times 10^{-2}$ & $3.17 \times 10^{-3}$ & $8.87 \times 10^{-7}$ &  $9.63 \times 10^{-9}$ \\
$R = 8$ & $1.14 \times 10^{-1}$ & $1.86 \times 10^{-2}$ & $6.31 \times 10^{-5}$ & $1.67 \times 10^{-11}$ \\
$R = 10$ & $5.56 \times 10^{-1}$ & $5.82 \times 10^{-2}$ & $9.94\times 10^{-4}$ & $1.12 \times 10^{-8}$
\end{tabular}}
\end{table}

\begin{table}[!ht]
\centering
\caption{$L^2$ difference between two methods computing the linearized collision term (Case 2)}
\label{tab:L2diff2}
\scalebox{.9}{%
\begin{tabular}{c|cccc}
& $N = 4$ & $N = 8$ & $N = 16$ & $N = 32$ \\
\hline
$R = 4$ & $1.87 \times 10^{-2}$ & $2.47 \times 10^{-4}$ & $2.30 \times 10^{-4}$ & $2.31 \times 10^{-4}$\\
$R = 6$ & $3.79 \times 10^{-1}$ & $1.92 \times 10^{-3}$ & $3.96 \times 10^{-5}$ & $2.17 \times 10^{-6}$ \\
$R = 8$ & $2.16$ & $2.39 \times 10^{-2}$ & $1.38 \times 10^{-4}$ & $6.63 \times 10^{-6}$\\
$R = 10$ & $7.34$ & $1.32 \times 10^{-1}$ & $5.16\times 10^{-4}$ & $1.50 \times 10^{-5}$
\end{tabular}}
\end{table}

We would like to point out that such a comparison only validates our method by measuring the similarity between our results and the results of a previous method. The numbers in the tables do not indicate the actual numerical error. For instance, Table \ref{tab:L2diff2} shows the $L^2$ difference between two methods is $2.30 \times 10^{-4}$ for $R = 4$ and $N = 16$. However, if we compare the linearized collision terms computed using our approach with $R = 4$, $N = 16$ with the one computed $R = 8$, $N = 32$, the $L^2$ difference is $3.76 \times 10^{-3}$, more than ten times larger than $2.30 \times 10^{-4}$. Here the result for $R = 8$ and $N = 32$ can be considered as a reference solution due to its larger value of $R$, meaning that the actual numerical error for the result of $R = 4$ and $N = 16$ probably has the magnitude $10^{-3}$ or even larger, due to the inadequate size of the velocity domain and the discontinuity in $f^{(2)}$. Nevertheless, the spectral accuracy of our method is still confirmed due to the fast decaying difference between our method and an existing method that has already been confirmed to possess the spectral accuracy.

Now we test the closeness between the linearized collision operator and the binary collision operator by computing the relative difference $\|\mathcal{L}[f] - \mathcal{Q}[f,f]\| / \|\mathcal{Q}[f,f]\|$ in both Case 1 and Case 2 under various parameter settings.
The binary collision operator $\mathcal{Q}[f,f]$ is again computed using the method introduced in \Cref{sec:fsm_binary}.
The results are presented in \Cref{tab:diff_lin_bin}.
In both cases, when $R$ is too small or $R/N$ is too large, the computed collision terms are unreliable.
For parameters that can yield accurate computations, one can observe that the difference between two collision terms is around 15\% in Case 1 and 7\% in Case 2.
Note that the two distribution functions in the test are not close to the Maxwellian, so we do not expect a very small relative difference in this table.
In the next subsection, we will demonstrate that these two collision operators will generate close results in the spatially homogeneous Boltzmann equation.

\begin{table}[!ht]
\centering
\caption{The relative $L^2$ difference between the linearized collision term and the binary collision term}
\label{tab:diff_lin_bin}
\scalebox{.9}{%
\begin{tabular}{c|cccc|cccc}
\hline
& \multicolumn{4}{|c}{Case 1}
& \multicolumn{4}{|c}{Case 2} \\
\cline{2-9}
& $N = 4$ & $N = 8$ & $N = 16$ & $N = 32$
& $N = 4$ & $N = 8$ & $N = 16$ & $N = 32$ \\
\hline
$R = 4$ & $0.287$ & $0.158$ & $0.159$ & $0.159$ & $1.030$ & $0.111$ & $0.090$ & $0.082$ \\
$R = 6$ & $0.643$ & $0.248$ & $0.150$ & $0.150$ & $3.73$ & $0.114$ & $0.080$ &  $0.070$ \\
$R = 8$ & $0.452$ & $0.384$ & $0.150$ & $0.150$ & $10.247$ & $1.097$ & $0.089$ & $0.073$ \\
$R = 10$ & $6.300$ & $0.592$ & $0.156$ & $0.150$ & $54.482$ & $3.476$ & $0.100$ & $0.076$ \\
\hline
\end{tabular}}
\end{table}

\subsection{Spatially homogeneous Boltzmann equation}
We now test the performance of our algorithm by applying the method to the spatially homogeneous Boltzmann equation. We will solve both equations with the quadratic and linearized operators:
\begin{equation} \label{eq:spatially_homogeneous_equations}
\frac{\partial f_{\mathrm{bin}}}{\partial t} = \mQ[f_{\mathrm{bin}},f_{\mathrm{bin}}], \qquad \frac{\partial f_{\mathrm{lin}}}{\partial t} = \mL[f_{\mathrm{lin}}],
\end{equation}
and compare both numerical results and the computational time. We will again consider two cases:
\begin{description}
    \item[Case 1] Collision kernel: $\mB_1$, initial distribution function: $f^{(1)}$;
    \item[Case 2] Collision kernel: $\mB_2$, initial distribution function $f^{(2)}$.
\end{description}
In both cases, the Maxwellian $\mM$ will be chosen as the standard normal distribution in the linearized collision operator. 

\begin{figure}[!ht]
\centering
\subfigure[Case 1]{
\includegraphics[scale=.35]{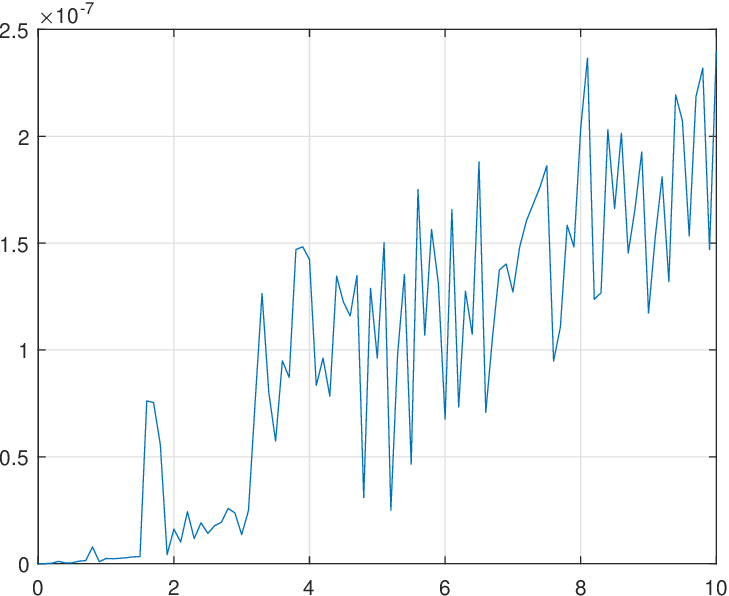}
} \qquad
\subfigure[Case 2]{
\includegraphics[scale=.35]{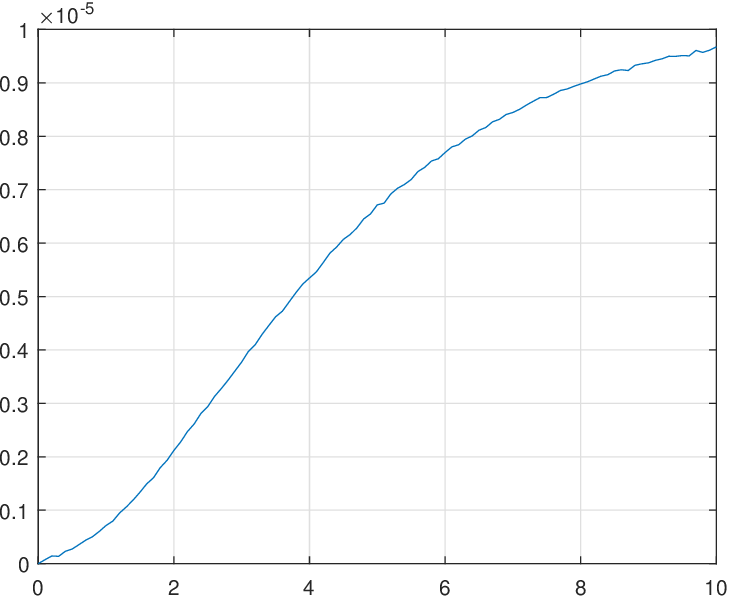}
}
\caption{$L^2$ difference between $f_{\mathrm{lin}}$ computed using double-precision floating-point numbers and quadruple-precision floating-point numbers.}
\label{fig:double_vs_binary128}
\end{figure}

In our implementation, the time derivative is discretized using the 4th-order Runge-Kutta scheme with time step $0.1$. In the spectral methods for the collision terms, the integral with respect to $g$ is again approximated using the Jacobi-Gauss-Radau quadrature that is exact for all polynomials up to degree $2N$. For the binary collision operator, the quadrature on the sphere is chosen with a lower accuracy to be specified later. In all the tests below, we will choose $R = 6$ and $N = 16$, which gives reasonable numerical accuracy according to the tests in the previous section. The purpose of the experiments is to show that
\begin{enumerate}
\item The linearized collision operator provides results similar to the results of the binary collision operator;
\item The computation of the linearized collision operator runs faster than the binary collision operator.
\end{enumerate}
Before comparing the two operators, we claim that the computation will fail in the very first time step if double-precision floating-point numbers are used without the cutoff \eqref{eq:cutoff}. The other two methods (using \texttt{double} with cutoff \eqref{eq:cutoff} \& using \texttt{binary128} without cutoff) have been both tested, and they have negligible differences in the numerical solutions. In Figure \ref{fig:double_vs_binary128}, we show the $L^2$ difference between the results computed with these two datatypes. In Case 1, the initial distribution function has a fast decay, so that even if the distribution function is divided by the Maxwellian, no large values are produced. As a result, the $L^2$ difference between two results is generally random. However, in Case 2, $f^{(2)}(\bv) / \mM(\bv)$ increases to infinity when $|\bv|$ tends to infinity in certain directions, so that the truncation has a more significant effect. Nevertheless, the difference is below $10^{-5}$ and will not affect the general conclusion below. Therefore, in the following tests, all results given are based on computations using double-precision floating-point numbers with cutoff \eqref{eq:cutoff}, unless otherwise specified.

\begin{figure}[!ht]
\centering
\subfigure[$t = 1$]{
\includegraphics[bb=50 20 380 300,clip,width=.3\textwidth]{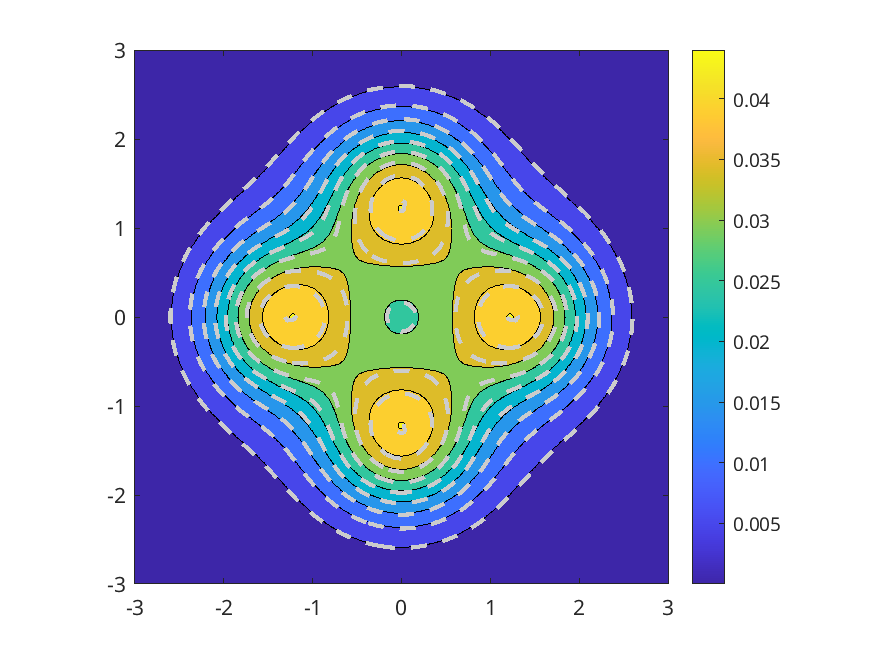}
}
\subfigure[$t = 2$]{
\includegraphics[bb=50 20 380 300,clip,width=.3\textwidth]{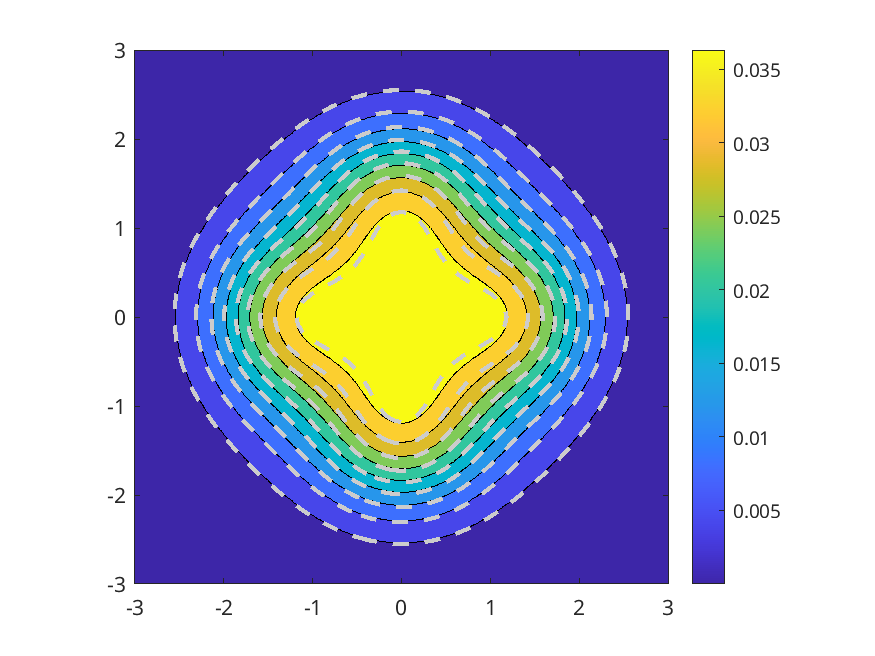}
}
\subfigure[$t = 4$]{
\includegraphics[bb=50 20 380 300,clip,width=.3\textwidth]{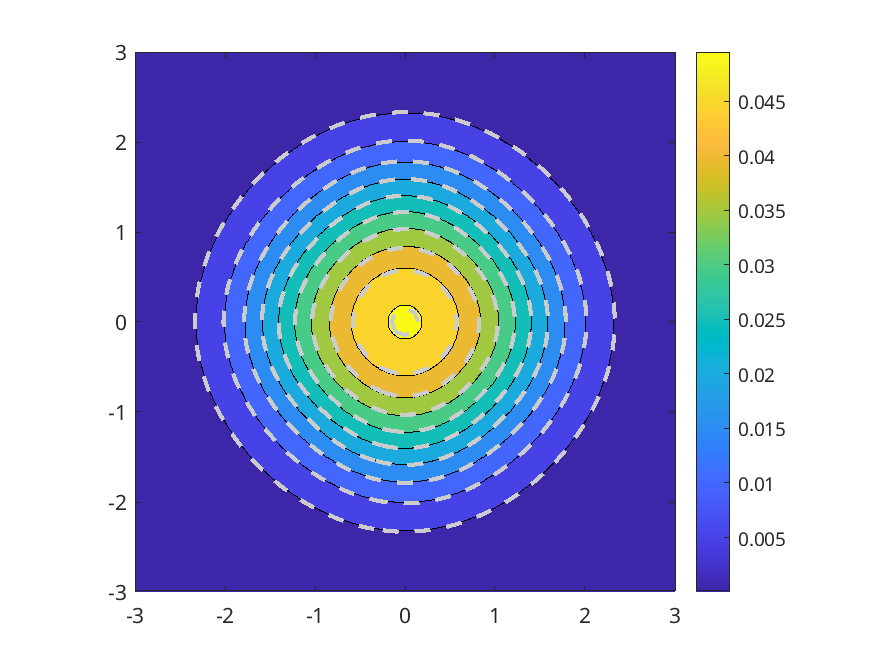}
}
\caption{Comparison between the numerical results of the binary and linearized collision operators for Case 1 at $v_3 = 0$. The filled contour plots with black contour lines represent the solutions of the linearized Boltzmann equation, and the gray dashed contour lines denote the solutions of the original Boltzmann equation.}
\label{fig:bin_vs_lin_Case1}
\end{figure}

\begin{figure}[!ht]
\centering
\subfigure[$t = 0.5$]{
\includegraphics[bb=50 20 380 300,clip,width=.3\textwidth]{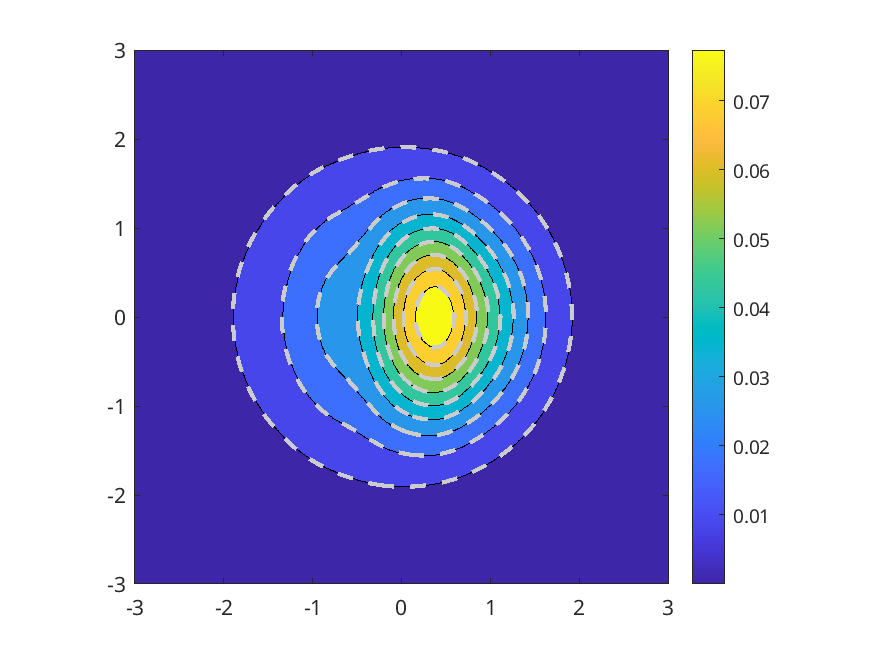}
}
\subfigure[$t = 2$]{
\includegraphics[bb=50 20 380 300,clip,width=.3\textwidth]{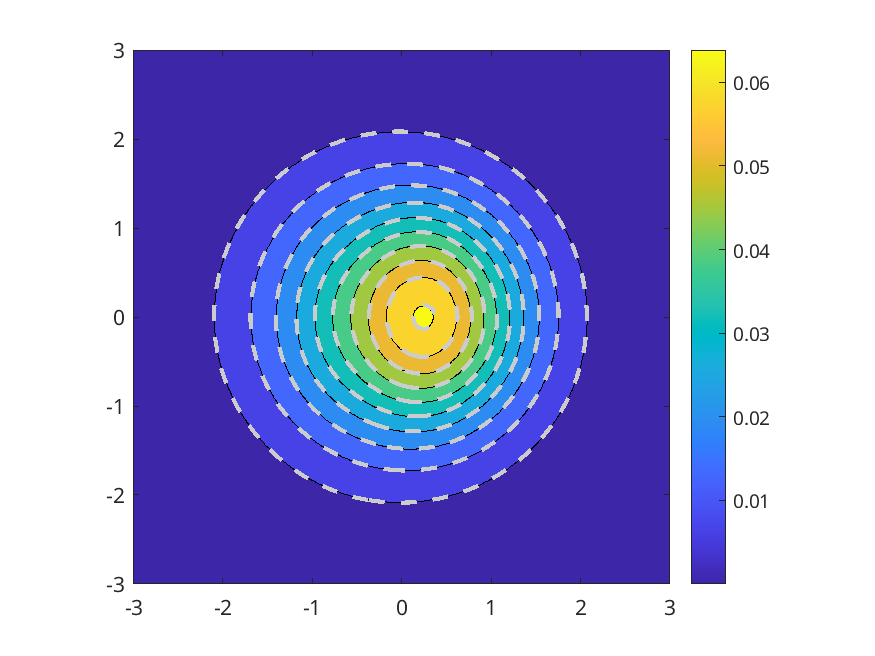}
}
\subfigure[$t = 6$]{
\includegraphics[bb=50 20 380 300,clip,width=.3\textwidth]{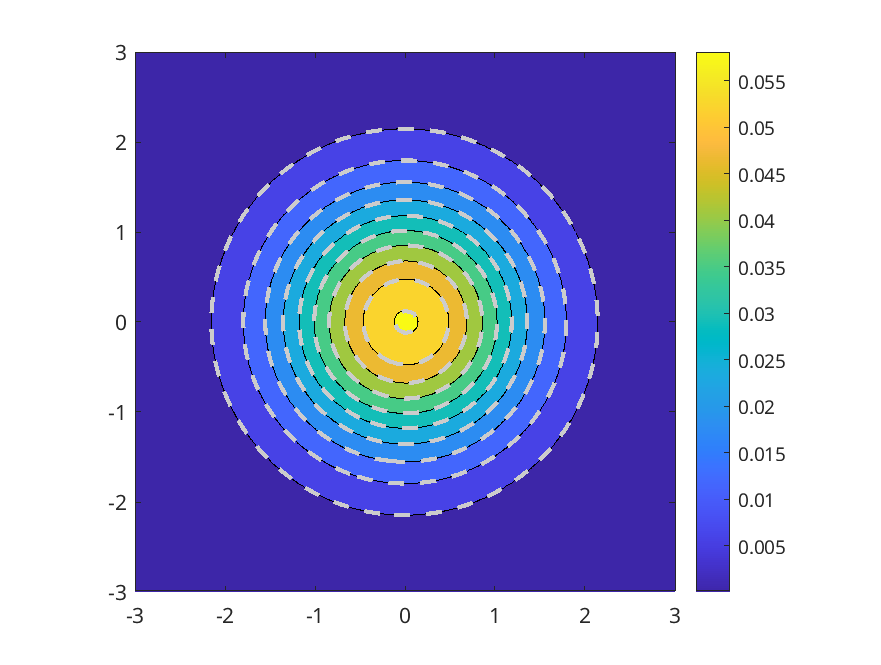}
}
\caption{Comparison between the numerical results of the binary and linearized collision operators for Case 2 at $v_3 = 0$. The filled contour plots with black contour lines represent the solutions of the linearized Boltzmann equation, and the gray dashed contour lines denote the solutions of the original Boltzmann equation.}
\label{fig:bin_vs_lin_Case2}
\end{figure}

We first carry out tests for implementations of the collision operators without conservation fixes. The contour lines of solutions to both equations at various times are given in Figures \ref{fig:bin_vs_lin_Case1} and \ref{fig:bin_vs_lin_Case2}. Note that the solutions are a three-dimensional functions. For better visualization, we only plot the contours of the slices at $v_3 = 0$. The numerical solution of $f_{\mathrm{bin}}$ is obtained by using 25 quadrature points on the hemisphere, and the linearized equation is solved with double-precision floating-point numbers and a cutoff. Here we remark that there is no visible difference between the contour lines of the $f_{\mathrm{lin}}$ obtained using both datatypes. It can be seen that the two sets of contour lines coincide in most of the figures, showing the similarity between two collision operators. Both solutions show an isotropic structure for large $t$, indicating the convergence towards the Maxwellian distribution function.

\begin{figure}[!ht]
\centering
\subfigure[Case 1]{
\includegraphics[scale=.35]{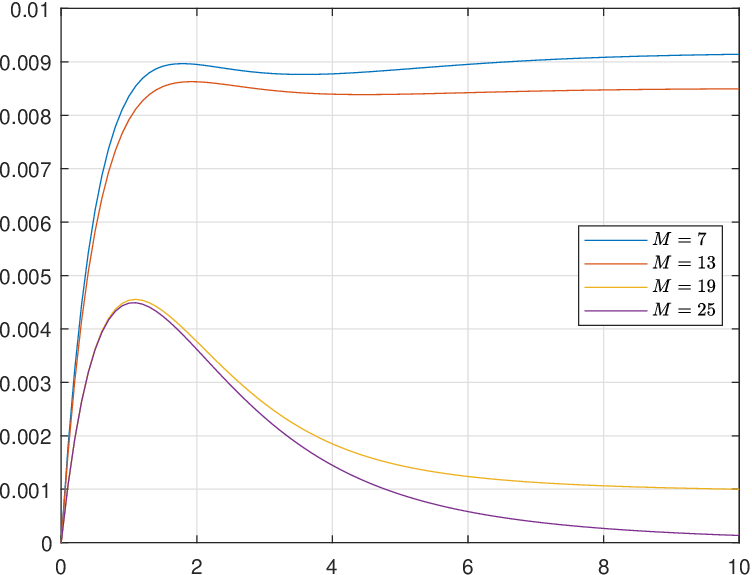}
} \qquad
\subfigure[Case 2]{
\includegraphics[scale=.35]{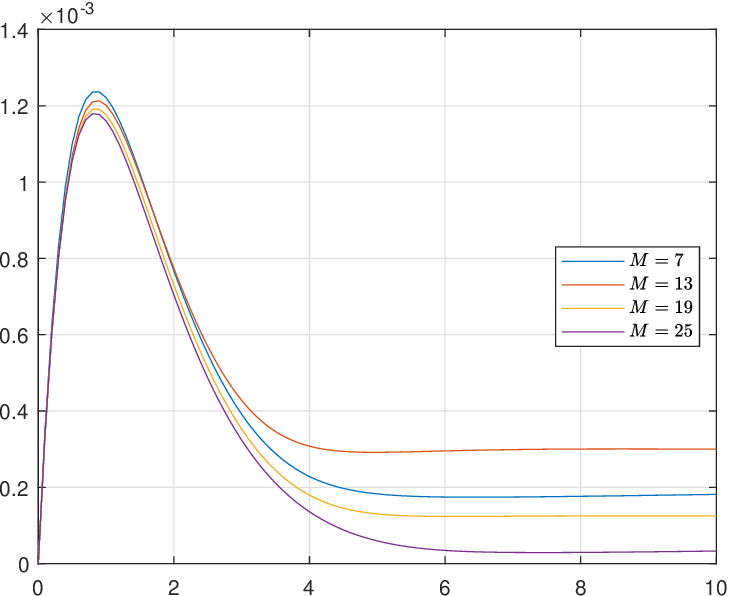}
}
\caption{$L^2$ difference between the results of linearized and binary collision operators. $M$: the number of quadrature points on the sphere.}
\label{fig:L2diff}
\end{figure}

\begin{figure}[!ht]
\centering
\includegraphics[scale=.35]{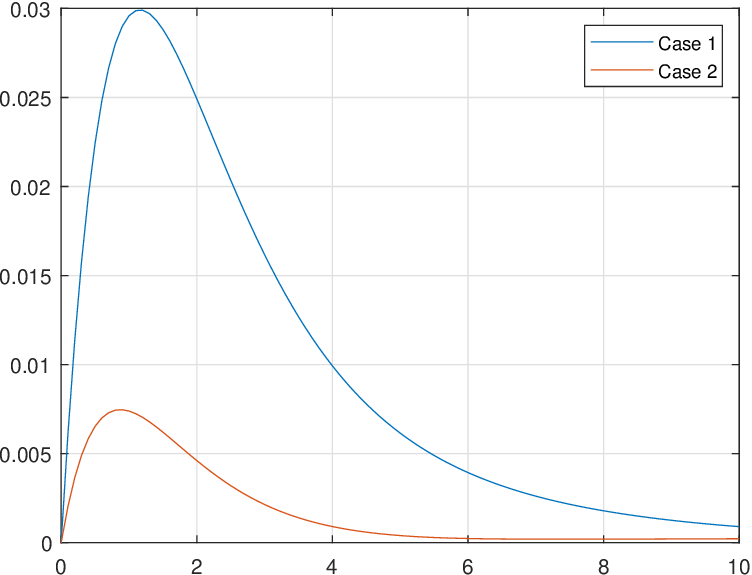}
\caption{The relative $L^2$ difference defined in \eqref{eq:Et}, where the result of the binary collision operator is computed using $M = 25$.}
\label{fig:RelDiff}
\end{figure}

For a more quantitative comparison, we plot the evolution of the $L^2$ difference between the solutions of two equations in \eqref{eq:spatially_homogeneous_equations}. Since the results for the binary collision operators depend on the choice of the numerical integration on the sphere, we have included four different quadrature rules with the number of quadrature points $M$ being $7$, $13$, $19$ and $25$. The difference between $f_{\mathrm{lin}}$ and $f_{\mathrm{bin}}$ is plotted in Figure \ref{fig:L2diff}. The results show that for the most accurate result of the binary collision term, given by $M = 25$, the $L^2$ difference between the solutions of the two equations is less than $0.005$ for Case 1 and less than $0.0013$ for Case 2. The relative $L^2$ difference defined by
\begin{equation} \label{eq:Et}
E(t) = \frac{\|f_{\mathrm{lin}}(t) - f_{\mathrm{bin}}(t)\|_{L^2}}{\|f_{\mathrm{bin}}(t)\|_{L^2}}
\end{equation}
is plotted in Figure \ref{fig:RelDiff}, showing that the difference never exceeds $3\%$. Figure \ref{fig:L2diff} also indicates that a sufficient large $M$ is required to get accurate results for the binary collision term.

Another way to confirm the correctness of the simulation is to check whether the distribution function indeed evolves towards the Maxwellian. In Table \ref{tab:l2diff_max}, we tabulate the $L^2$ difference between the Maxwellian and the solution at $t = 10$. We can again observe that the accuracy of our method for the linearized collision operator is comparable to the result for the binary collision term with $M = 25$. Further reduction of the difference requires larger values of $R$ and $N$.
\begin{table}[!ht]
\centering
\caption{$L^2$ difference between the Maxwellian and the solution at $t = 20$.}
\label{tab:l2diff_max}
\scalebox{.9}{%
\begin{tabular}{c@{\qquad}ccccc}
\hline
& $M = 7$ & $M = 13$ & $M = 19$ & $M = 25$ & Linearized \\
\hline
Case 1 & $9.15 \times 10^{-3}$ & $8.48 \times 10^{-3}$ & $9.60 \times 10^{-4}$ & $4.90 \times 10^{-5}$ & $3.09 \times 10^{-5}$ \\
Case 2 & $4.44 \times 10^{-4}$ & $6.43 \times 10^{-4}$ & $4.81 \times 10^{-4}$ & $4.25 \times 10^{-4}$ & $4.25 \times 10^{-4}$ \\
\hline
\end{tabular}
}
\end{table}

\begin{figure}[!ht]
\centering
\subfigure[Case 1]{
\includegraphics[scale=.35]{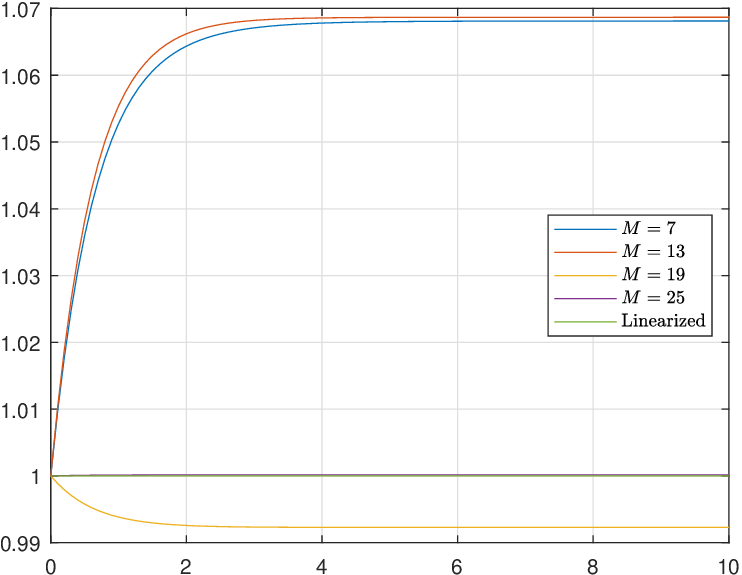}
} \quad
\subfigure[Case 2]{
\includegraphics[scale=.35]{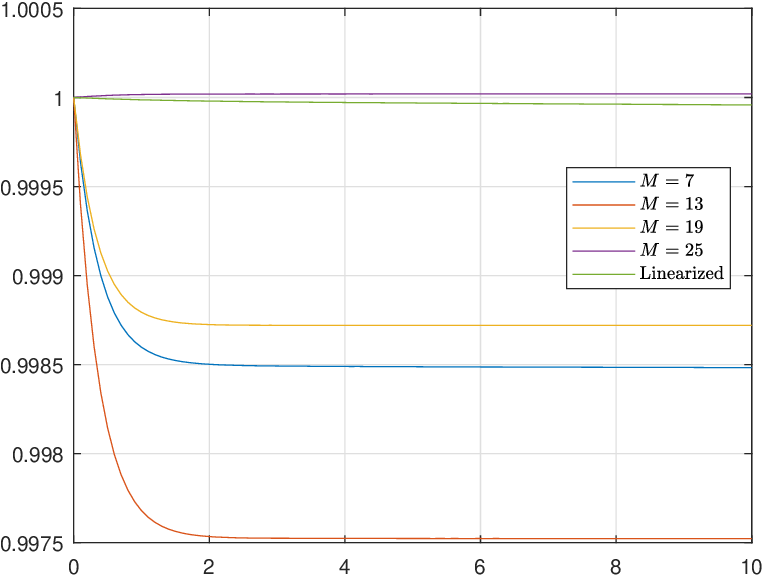}
}
\caption{Evolution of the mass for both equations in \eqref{eq:spatially_homogeneous_equations}.}
\label{fig:mass}
\end{figure}

\begin{figure}[!ht]
\centering
\subfigure[Case 1]{
\includegraphics[scale=.35]{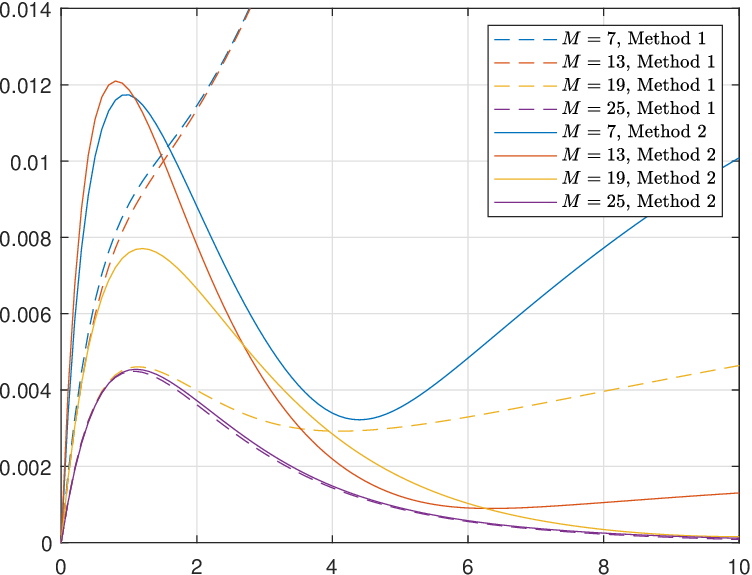}
} \quad
\subfigure[Case 2]{
\includegraphics[scale=.35]{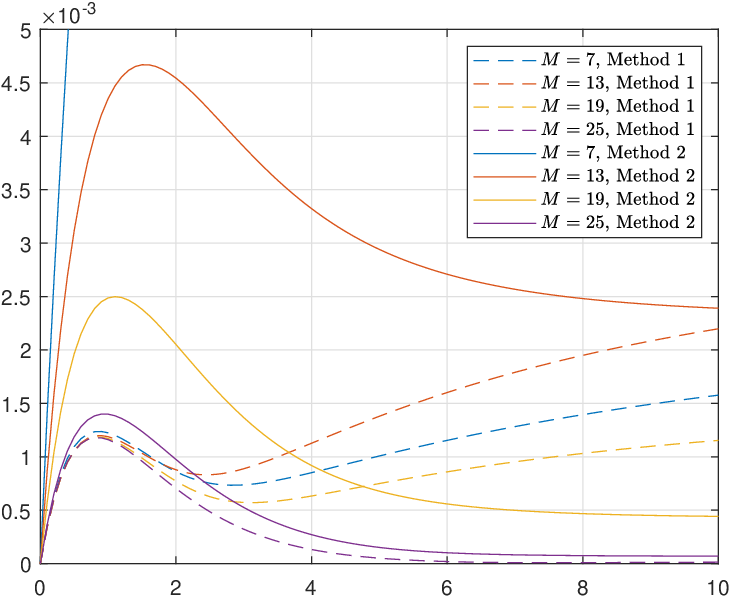}
}
\caption{$L^2$ difference between the results of linearized and binary collision operators. For the linearized equation, the conservation of mass is enforced by setting $\widehat{\mL}_0 = 0$ at every time step. For the quadratic equation, this is enforced by two different methods. Method 1 refers to setting $\hat{Q}_0$ to $0$, and Method 2 refers to replacing the sinc function with \eqref{eq:sinc}.}
\label{fig:l2_diff_csv}
\end{figure}

Another phenomenon we have observed in the experiments is the loss of mass conservation in the numerical results of both equations. In Figure \ref{fig:mass}, we plot the evolution of the masses in all test cases. In general, for the linearized collision operator, our method does not maintain exact mass conservation, but the error is below 0.01\% in both cases at $t = 10$, which is acceptable in most applications. In the binary case, mass conservation relies on an accurate integration on the sphere (see \eqref{eq:hatQ0}), and therefore a larger $M$ gives better mass conservation. Such an observation inspires us to check if enforcing mass conservation using methods introduced at the end of Section \ref{sec:fsm_binary} can provide better accuracy, especially for smaller values of $M$ in the computation of binary collision operators. We therefore apply the fix of mass conservation to the numerical solvers of both equations. For the linearized equation, the mass conservation is enforced by setting $\widehat{\mL}_0 = 0$ every time the collision operator is evaluated. For the quadratic equation, the two methods discussed at the end of Section \ref{sec:fsm_binary} are tested. We again study the $L^2$ difference between numerical solutions of the two equations in \eqref{eq:spatially_homogeneous_equations}, and the results are given in Figure \ref{fig:l2_diff_csv}. For the linearized equation and the quadratic equation solved with $M = 25$, the conservation fix does not have much effect on the numerical solution, since the original method already preserves the total mass well. However, comparing Figure \ref{fig:L2diff} and Figure \ref{fig:l2_diff_csv}, we find that imposing mass conservation does not help reduce the overall error of the solution, especially for long-time simulations. In particular, for $M = 7$ in Case 1, the solution quickly loses validity with both fixing methods and does not appear to approach the correct equilibrium state. This indicates the necessity of a sufficiently large $M$ in some simulations.

We now compare the computational time of different numerical methods. Since the conservation fix only takes a tiny proportion of the total computational time, below we only provide the time used for methods without the conservation fix. The tests are carried out on the Intel\textregistered{} Core\texttrademark{} i7-13800H processor without parallelization, and the computational time for 100 Runge-Kutta steps is provided in Table \ref{tab:computational_time}. In general, all operations for Case 1 and Case 2 are the same, and the differences in computational times in both rows are mainly due to external reasons such as background processes. As expected, for the quadratic equation, the computational time is roughly proportional to the value of $M$, and the time for the linearized equation with datatype \texttt{double} is slightly longer than the time for the quadratic equation divided by $M$. However, the datatype \texttt{binary128} is significantly slower on the processor tested, which is not worth considering in practical uses. 
\begin{table}[!ht]
\centering
\caption{Computational time for solving the spatially homogeneous equation from $t = 0$ to $t = 10$ with time step $0.1$.}
\label{tab:computational_time}
\scalebox{.9}{%
\begin{tabular}{c@{\qquad}ccccccc}
\hline
& \multicolumn{4}{c}{Quadratic equation} & & \multicolumn{2}{c}{Linearized equation} \\
\cline{2-5} \cline{7-8}
& $M = 7$ & $M = 13$ & $M = 19$ & $M = 25$ & & \texttt{double} & \texttt{binary128} \\
\hline
Case 1 & 15.29s & 27.58s & 40.02s & 52.53s & & 2.81s & 156.14s \\
Case 2 & 14.75s & 27.17s & 40.21s & 52.25s & & 2.81s & 155.22s \\
\hline
\end{tabular}}
\end{table}

\section{Application to the steady-state Boltzmann equation}
\label{sec:application}
The most direct application of the linearized Boltzmann collision operator is to use it to replace the binary collision operator and accelerate simulations. This direction will be explored in our future work, and here we will use our algorithm to help accelerate the numerical solver of the steady-state Boltzmann equation with the quadratic collision term.

\subsection{Problem settings}
For simplicity, we focus on the steady-state Boltzmann equation with a one-dimensional spatial variable $x$ and a three-dimensional velocity variable $\bv$. The equation reads
\begin{equation} \label{eq:steady_Boltz}
v_1 \partial_x f = \frac{1}{\Kn} \mQ[f,f], \qquad x \in (x_L, x_R), \quad \bv \in \mathbb{R}^3.
\end{equation}
Here, $v_1$ refers to the first component of $\bv$, and we have added a dimensionless parameter $\Kn$ in front of the collision term to indicate the rarefaction of the gas. The collision kernel is chosen to be $\mB_1$ (see \eqref{eq:col_ker}), so that the gas molecules are Maxwell molecules. Maxwell's wall boundary conditions with full accommodation are imposed at $x_L$ and $x_R$:
\begin{equation} \label{eq:maxwell_bc}
\begin{aligned}
f(x_L, \bv) &= \frac{\rho_L}{(2\pi \theta_L)^{3/2}} \exp \left( -\frac{|\bv - \bu_L|^2}{2 \theta_L} \right), \qquad \text{if } v_1 > 0, \\
f(x_R, \bv) &= \frac{\rho_R}{(2\pi \theta_R)^{3/2}} \exp \left( -\frac{|\bv - \bu_R|^2}{2 \theta_R} \right), \qquad \text{if } v_1 < 0,
\end{aligned}
\end{equation}
where we have used $\bu_L$ and $\theta_L$ to denote the velocity and temperature of the left wall, and $\bu_R$ and $\theta_R$ correspond to the right wall. We assume that the first components of both $\bu_L$ and $\bu_R$ are zero, so that the domain $(x_L, x_R)$ stays unchanged. The ``densities'' $\rho_L$ and $\rho_R$ are chosen to guarantee there is no mass flux on the wall. The explicit expresssions are
\begin{displaymath}
\rho_L = -\sqrt{\frac{2\pi}{\theta_L}} \int_{\mathbb{R}^3} \min(v_1,0) f(x_L, \bv) \,\dd \bv, \quad
\rho_R = \sqrt{\frac{2\pi}{\theta_R}} \int_{\mathbb{R}^3} \max(v_1,0) f(x_L, \bv) \,\dd \bv.
\end{displaymath}
To uniquely determine the solution, an additional condition specifying the total mass is needed:
\begin{equation}
\label{eq:total_mass}
\int_{x_L}^{x_R} \int_{\mathbb{R}^3} f(x,\bv)  \,\dd\bv \,\dd x = C,
\end{equation}
where $C$ is a given constant.

\subsection{Modified Newton's method}
We are going to develop an iterative solver for the steady-state Boltzmann equation \eqref{eq:steady_Boltz}. To achieve fast convergence, we consider Newton's iteration defined by
\begin{equation}
\label{eq:Newton_iter}
f^{(k+1)} = f^{(k)} - g^{(k)},
\end{equation}
where $g^{(k)}$ satisfies the linear equation
\begin{equation}
\label{eq:lin_eq}
v_1 \partial_x (g^{(k)} - f^{(k)}) = \frac{1}{\Kn} (\mQ[f^{(k)}, g^{(k)}] + \mQ[g^{(k)}, f^{(k)}] - \mQ[f^{(k)}, f^{(k)}]).
\end{equation}
Since the boundary condition of $f$ (see \eqref{eq:maxwell_bc}) is linear, the function $g^{(k)}$ should satisfy the same boundary condition as that of $f$. Although Newton's method has a locally quadratic convergence rate, solving the linear equation \eqref{eq:lin_eq} can be challenging due to the large computational cost of $\mQ[f^{(k)}, g^{(k)}] + \mQ[g^{(k)}, f^{(k)}]$. In general, solving \eqref{eq:maxwell_bc} requires an iterative scheme, which we will refer to as ``inner iteration'' nested inside the Newton iteration. For most iterative methods for linear equations, each inner iteration requires at least one evaluation of $\mQ[f^{(k)}, g^{(k)}] + \mQ[g^{(k)}, f^{(k)}]$, which will harm the high efficiency of Newton's method.

Note that although $\mQ[f^{(k)}, g^{(k)}] + \mQ[g^{(k)}, f^{(k)}]$ is also a linearization of the collision operator, our fast algorithm does not apply since the linearization is not about a local equilibrium state. In order to apply our fast algorithm, we replace the equation of $g^{(k)}$ \eqref{eq:lin_eq} with
\begin{equation}
\label{eq:modified_lin_eq}
v_1 \partial_x (g^{(k)} - f^{(k)}) = \frac{1}{\Kn} (\mL^{(k)}[g^{(k)}] - \mQ[f^{(k)}, f^{(k)}]),
\end{equation}
where $\mL^{(k)}$ is defined by
\begin{displaymath}
\mL^{(k)}[g^{(k)}] = \mQ[\mathcal{M}^{(k)}, g^{(k)}] + \mQ[g^{(k)}, \mathcal{M}^{(k)}],
\end{displaymath}
and $\mathcal{M}^{(k)}$ is the local Maxwellian associated with $f^{(k)}$. 
This approach is similar to the quasi-Newton method, which replaces the exact Jacobian with an approximate Jacobian, but still guarantees the exact solution is a fixed point of the iterative method. Thus, for each Newton iteration, the quadratic operator needs to be evaluated only once on each spatial point when computing $\mQ[f^{(k)}, f^{(k)}]$ in \eqref{eq:modified_lin_eq}.

The modified Newton's method \eqref{eq:Newton_iter}\eqref{eq:modified_lin_eq} still requires solving a linear system. In this work, we adopt a classical method known as ``source iteration'' to solve $g^{(k)}$ in \eqref{eq:modified_lin_eq}. Let $g^{(k,\ell)}$ be the numerical solution of $g^{(k)}$ in the $\ell$th iteration. With the initial condition $g^{(k,0)}$ being zero, the source iteration updates the solution according to
\begin{equation} \label{eq:source_iteration}
v_1 \partial_x (g^{(k,\ell+1)} - f^{(k)}) + \frac{\nu}{\Kn} g^{(k,\ell+1)} = \frac{1}{\Kn} (\mL^{(k)}[g^{(k,\ell)}] + \nu g^{(k,\ell)} - \mQ[f^{(k)}, f^{(k)}]).
\end{equation}
The parameter $\nu$ is chosen such that $\mL^{(k)}[g^{(k,\ell)}] + \nu g^{(k,\ell)}$ is positive. 
Using the upwind method to discretize the spatial derivative, one can solve $g^{(k,\ell+1)}$ via fast sweeping. For instance, for a velocity $\bv_j$ with $v_{j,1} > 0$, the upwind method of \eqref{eq:source_iteration} reads
\begin{equation}
\begin{aligned}
v_{j,1} \left( \frac{g_{i,j}^{(k,\ell+1)} - g_{i-1,j}^{(k,\ell+1)}}{\Delta x} - \frac{f_{i,j}^{(k)} - f_{i-1,j}^{(k)}}{\Delta x} \right) + \frac{\nu_i}{\Kn} g_{i,j}^{(k,\ell+1)} = \qquad \\
\frac{1}{\Kn} (\mL_j^{(k)}[g_i^{(k,\ell)}] + \nu_i g_{i,j}^{(k,\ell)} - \mQ_j[f_i^{(k)}, f_i^{(k)}]).
\end{aligned}
\end{equation}
Here the subscript $i$ is the index of the spatial grid cell, and $j$ is the velocity index.
This equation shows that $g_{i,j}^{(k,\ell+1)}$ can be fully determined $g_{i-1,j}^{(k,\ell+1)}$.
The boundary condition is handled by the ghost cell method, but the function in the ghost cell is determined by the solution in the previous iterative step $g^{(k,\ell+1)}$, which can be computed before solving $g^{(k,\ell+1)}$.
Thus, we can start from the leftmost grid cell and perform a left-to-right scan to find out $g_{i,j}^{(k,\ell+1)}$ for all $i$.
Similarly, if $v_{j,1} < 0$, we can perform a right-to-left scan to obtain the solution.

The iteration \eqref{eq:source_iteration} is known to have slow convergence when $\Kn$ is small, and the convergence rate can be significantly improved by using the general synthetic iterative method \cite{Su2020}. Since our major purpose is to demonstrate that fast convergence can be still achieved using the modified Newton's method \eqref{eq:Newton_iter}\eqref{eq:modified_lin_eq}, we will leave the exploration of better inner iterations to future works. 

A fully discrete numerical scheme requires discretization of $x$ and $\bv$. In general, the spatial variable is discretized with a first-order finite volume scheme, and the velocity variable is discretized using the Fourier spectral method described in previous sections. Since the modified Newton's method is generally formulated at the continuous level, the generalization to higher-order schemes such as discontinuous Galerkin and WENO methods is straightforward. More details about our spatial discretization can be found in supplementary materials.

\subsection{Numerical results}
In our experiments, we choose $x_L = -0.5$ and $x_R = 0.5$. The total mass $C$ is set to be $1.0$. The spatial domain is gridded into $200$ cells, and $N = 16$ ($32$ modes in each direction) is used in velocity discretization. The choice of $R$ (truncation of the relative velocity) will depend on the temperature of the walls. In general, larger values of $R$ are needed for higher temperatures. We stop the computation if the $L^2$ norm of the residual defined by
\begin{equation} \label{eq:res_r}
\sqrt{\int_{x_L}^{x_R} \int_{[-L,L]^3} \left| v_1 \partial_x f - \frac{1}{\Kn} \mQ[f,f] \right|^2 \,\mathrm{d}\bv \,\mathrm{d}x}
\end{equation}
is less than $10^{-5}$. For inner iterations solving the system \eqref{eq:modified_lin_eq}, we terminate the iteration when the $L^2$ norm of the residual is either less than $10^{-7}$ or less than $10^{-3}$ times the initial residual. Since the initial guess of the solution to the linear system is zero, the initial residual is the same as \eqref{eq:res_r}. According to different choices of boundary conditions, two cases will be studied in the following subsections.

\subsubsection{Couette flow}
\begin{figure}[!ht]
\centering
\subfigure[Velocity]{
\includegraphics[scale=.35]{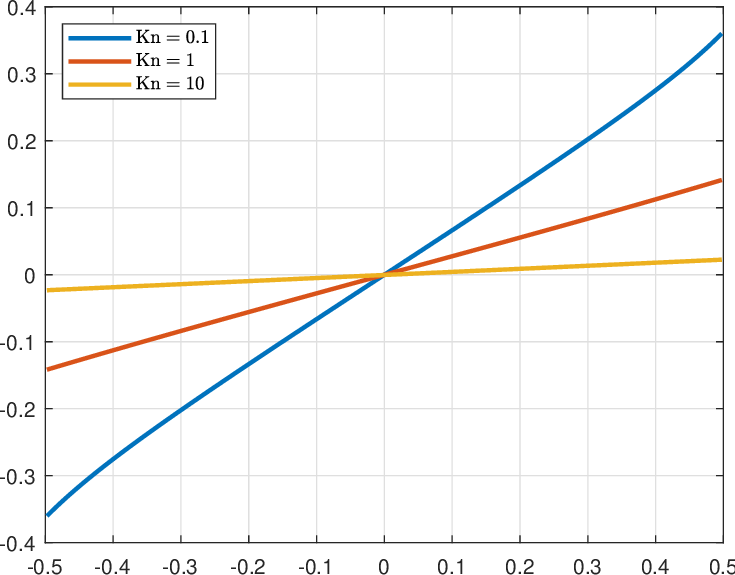}
} \quad
\subfigure[Temperature]{
\includegraphics[scale=.35]{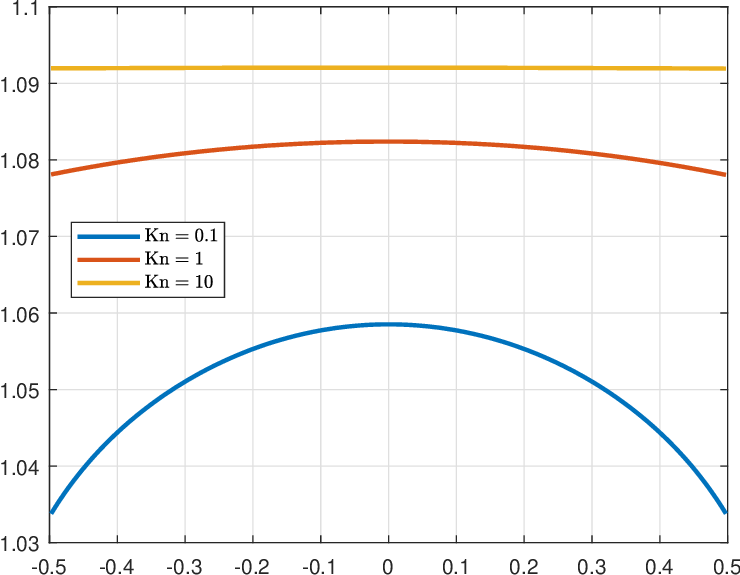}
}
\caption{Results of Couette flows with $u_W =0.5$.}
\label{fig:Couette_Kn}
\end{figure}

\begin{figure}[!ht]
\centering
\subfigure[Velocity]{
\includegraphics[scale=.35]{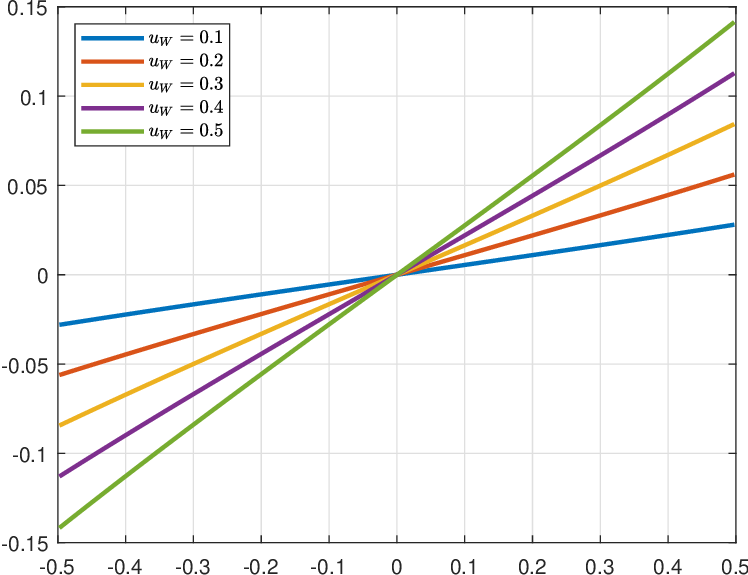}
} \quad
\subfigure[Temperature]{
\includegraphics[scale=.35]{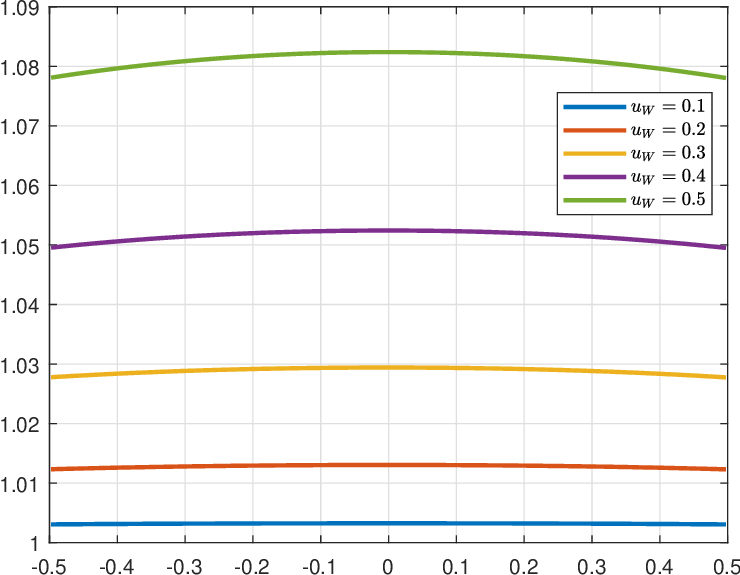}
}
\caption{Results of Couette flows with $Kn = 1.0$.}
\label{fig:Couette_uw}
\end{figure}

In our first example, we consider Couette flows where $\theta_L = \theta_R = 1$, and
\begin{displaymath}
\bu_L = (0, -u_W, 0)^T, \qquad \bu_R = (0, u_W, 0)^T.
\end{displaymath}
We first fix the wall speed $u_W$ to $0.5$ and study three different Knudsen numbers $\Kn = 0.1, 1$ and $10$. The parallel velocity and the temperature of the gases are plotted in Figure \ref{fig:Couette_Kn}. In general, for higher Knudsen numbers, both the velocity and temperature curves have smaller gradients, due to less interactions between the particles reflected from left and right walls. In particular, the steady-state distribution functions are independent of $x$ if $\Kn = \infty$. Figure \ref{fig:Couette_uw} shows the solution at Knudsen number $1.0$ for five different wall speeds. It is correctly predicted that higher wall speeds lead to larger velocity and higher temperature of the gas. 


\begin{figure}[!ht]
\centering
\subfigure[$Kn = 0.1$]{
\includegraphics[scale=.35]{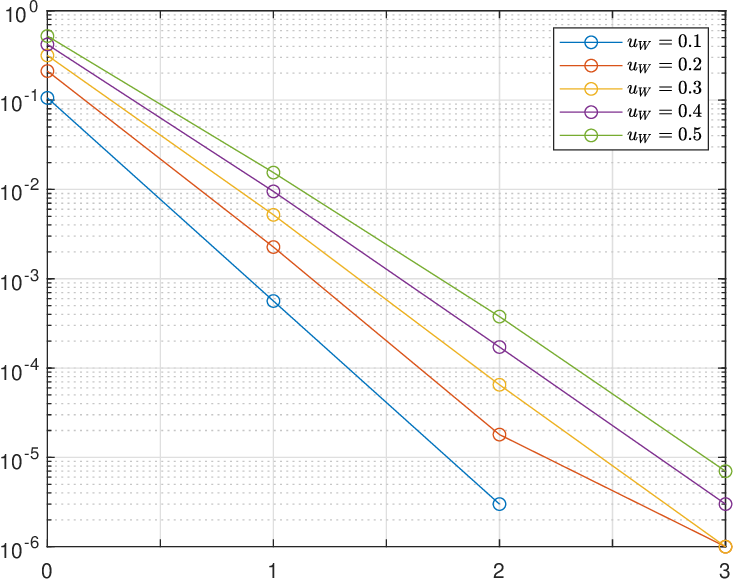}
} \quad
\subfigure[$Kn = 1$]{
\includegraphics[scale=.35]{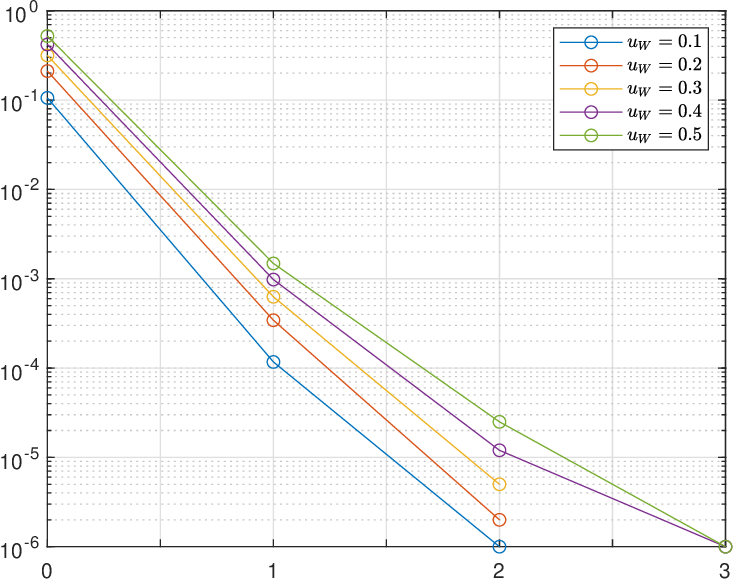}
}
\caption{Convergence of Newton iterations for Couette flows. The horizontal axis denotes the number of iterations, and the vertical axis stands for the $L^2$ residual.}
\label{fig:Couette_Newton}
\end{figure}

We now study the convergence behavior of our algorithm. In general, the modified Newton's method requires only a few iterations to converge. Table \ref{tab:Newton_cvg_Couette} shows the number of Newton iterations to achieve a residual with $L^2$ norm less than $10^{-5}$, which is also the number of evaluations of binary collision operators on each grid cell. For larger wall velocities and smaller Knudsen numbers, slightly more iterations are needed, since the solutions are generally further away from the initial value. Figure \ref{fig:Couette_Newton} provides the decay of the $L^2$ residual, which gives a clearer view of larger residuals at each iteration for larger wall velocities. The decay of the residual is faster for larger Knudsen numbers, due to the factor $1/\Kn$ in front of the collision term.
\begin{table}[!ht]
\centering
\caption{Numbers of Newton iterations for Couette flows}
\label{tab:Newton_cvg_Couette}
\scalebox{.9}{%
\begin{tabular}{c|ccccc}
  & $u_W = 0.1$ & $u_W = 0.2$ & $u_W = 0.3$ & $u_W = 0.4$ & $u_W = 0.5$ \\
  \hline
  $\Kn = 0.1$ & $2$ & $3$ & $3$ & $3$ & $3$ \\
  $\Kn = 1$ & $2$ & $2$ & $2$ & $3$ & $3$ \\
  $\Kn = 10$ & $2$ & $2$ & $2$ & $2$ & $2$
\end{tabular}
}
\end{table}

\begin{figure}[!ht]
\centering
\subfigure[$Kn = 0.1$]{
\includegraphics[scale=.35]{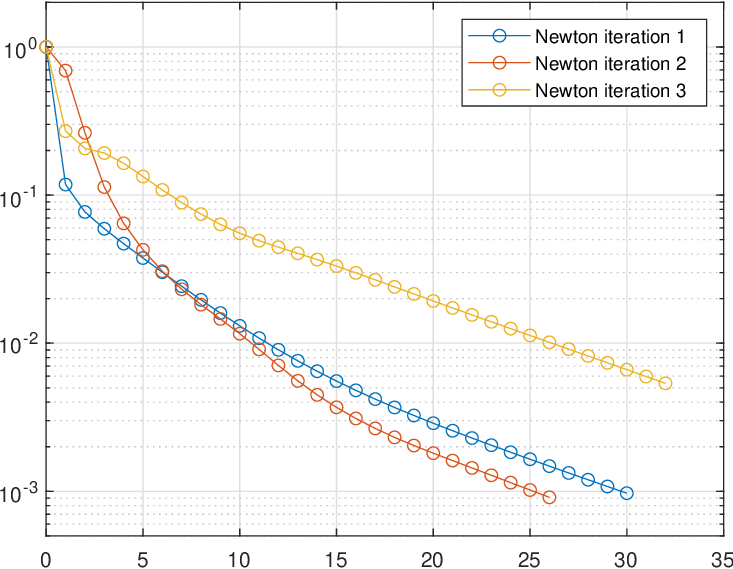}
} \quad
\subfigure[$Kn = 1$]{
\includegraphics[scale=.35]{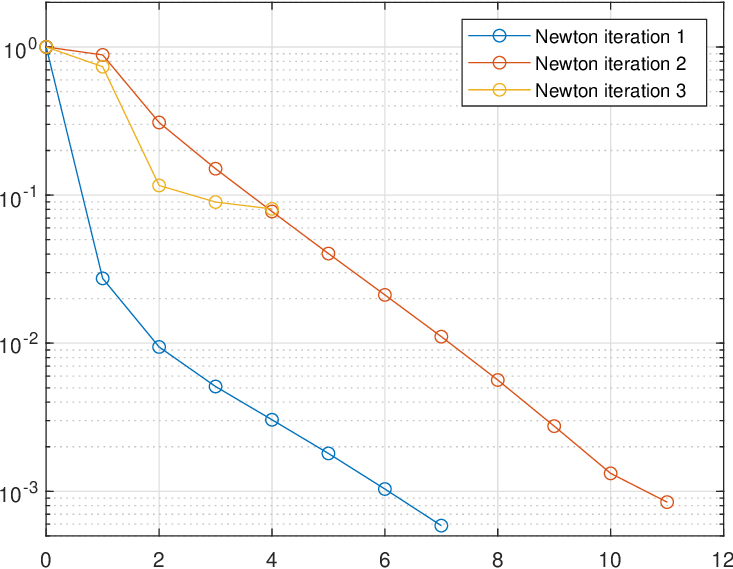}
}
\caption{Convergence of inner iterations for Couette flows with $u_W = 0.4$. The horizontal axis denotes the number of iterations, and the vertical axis stands for the relative $L^2$ residual (the residual at the current step divided by the initial residual).}
\label{fig:Couette_inner}
\end{figure}

For the source iteration, the convergence is significantly faster for larger Knudsen numbers. To show this, we take $u_W = 0.4$ as an example, for which both $\Kn = 0.1$ and $\Kn = 1$ require three Newton iterations. For all these Newton iterations, the decay of the residual is plotted in Figure \ref{fig:Couette_inner}. It is clearer that the case $\Kn = 0.1$ requires much more inner iterations to achieve the relative residual $10^{-3}$. Note that for the third Newton step in both figures, we terminate the inner iterations before the relative residual drops to $10^{-3}$, since the other threshold for the absolute residual, $10^{-7}$, has already been hit. As mentioned previously, improvements using micro-macro decomposition will be studied in future works. Nevertheless, each inner iteration only requires evaluation of one linearized collision operator on each grid cell, which is significantly faster than the evaluation of binary collision operators.
To support our arguments, we list in \Cref{tab:time} the computational times for calculating the residuals of the Boltzmann equation \eqref{eq:steady_Boltz} and the linear equation \eqref{eq:lin_eq}, based on the simulations run on the 11th Generation Intel\textsuperscript{\textregistered} Core\texttrademark{} i7 Processor with 8 threads.
The residual of \eqref{eq:steady_Boltz}, requiring evaluation of the binary collision operator, needs to be calculated only once in each Newton iteration, and the residual of \eqref{eq:lin_eq}, requiring the linearized collision operator, has to be evaluated for every inner iteration.
One can notice that even for $\Kn = 0.1$ that requires a large number of inner iterations, the computational time for the residual of \eqref{eq:steady_Boltz} still takes more than a half of the total computational time.
\begin{table}[!ht]
\centering
\caption{Computation times for evaluating the residuals in the simulations of Couette flows}
\label{tab:time}
\begin{tabular}{ccccc}
      \hline \hline
      && Time for the & Time for the & \\
      && residual of \eqref{eq:steady_Boltz} & residual of \eqref{eq:lin_eq} & Total \\
      \hline
      & $\Kn=0.1$  & 101.7\,s (58.2\%) & 60.0\,s (34.3\%) & 174.8\,s \\
      $u_W=0.1$ & $\Kn=1  $  & 107.6\,s (84.7\%) & 12.9\,s (10.1\%) & 127.1\,s \\
      & $\Kn=10 $  & 109.0\,s (89.6\%) & 6.8\,s  (5.6\%)  & 121.6\,s \\
      \hline
      & $\Kn=0.1$  & 147.5\,s (60.9\%) & 75.3\,s (31.1\%) & 242.3\,s \\
      $u_W=0.2$ & $\Kn=1  $  & 108.3\,s (80.8\%) & 14.6\,s (10.9\%) & 134.0\,s \\
      & $\Kn=10 $  & 108.5\,s (87.0\%) & 9.9\,s  (7.9\%)  & 124.7\,s \\
      \hline
      & $\Kn=0.1$  & 143.7\,s (56.4\%) & 84.2\,s (33.0\%) & 254.9\,s \\
      $u_W=0.3$ & $\Kn=1  $  & 107.2\,s (78.5\%) & 16.9\,s (12.4\%) & 136.5\,s \\
      & $\Kn=10 $  & 113.7\,s (84.3\%) & 15.0\,s (11.1\%) & 134.9\,s \\
      \hline
      & $\Kn=0.1$  & 142.4\,s (55.9\%) & 86.1\,s (33.8\%) & 254.9\,s \\
      $u_W=0.4$ & $\Kn=1  $  & 142.4\,s (78.7\%) & 21.3\,s (11.8\%) & 180.9\,s \\
      & $\Kn=10 $  & 110.8\,s (81.1\%) & 20.6\,s (15.1\%) & 136.7\,s \\
      \hline
      & $\Kn=0.1$  & 148.0\,s (57.7\%) & 85.5\,s (33.3\%) & 256.5\,s \\
      $u_W=0.5$ & $\Kn=1  $  & 148.8\,s (74.5\%) & 38.1\,s (19.1\%) & 199.6\,s \\
      & $\Kn=10 $  & 107.5\,s (70.2\%) & 32.9\,s (21.5\%) & 153.1\,s \\
      \hline \hline
\end{tabular}
\end{table}

\subsubsection{Fourier flow}
\begin{figure}[!ht]
\centering
\subfigure[Density]{
\includegraphics[scale=.35]{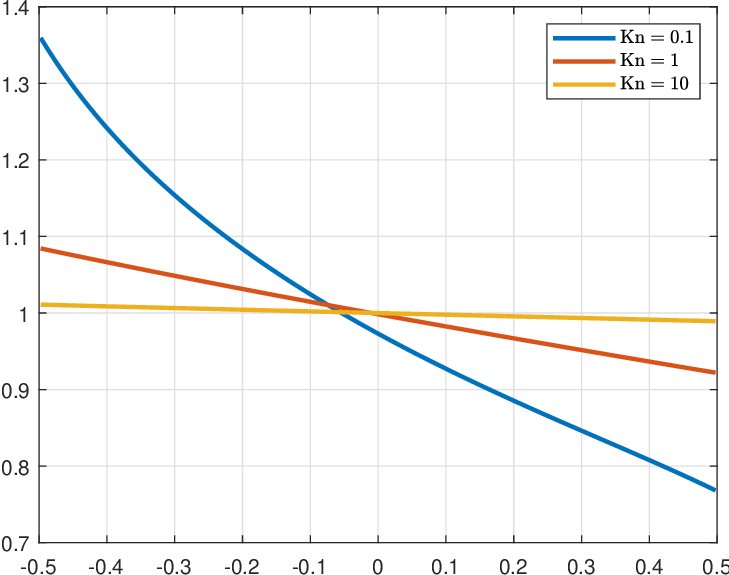}
} \quad
\subfigure[Temperature]{
\includegraphics[scale=.35]{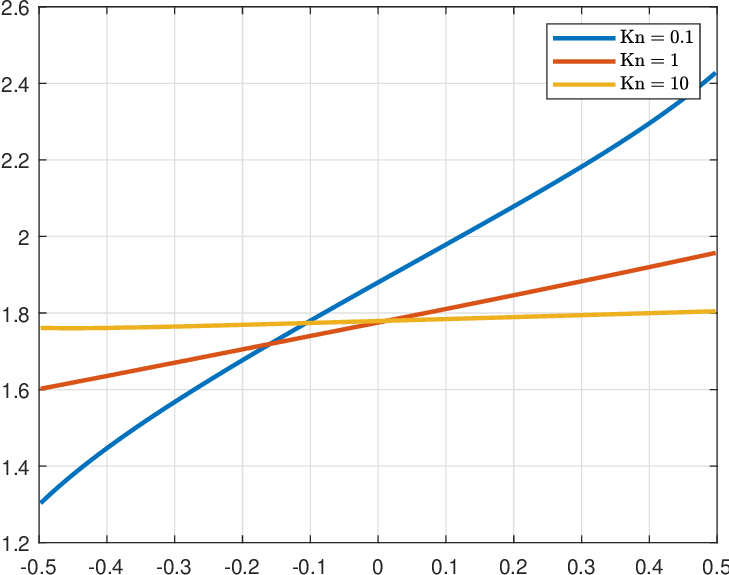}
}
\caption{Results of Fourier flows with $\theta_R = 3$.}
\label{fig:Fourier_Kn}
\end{figure}

\begin{figure}[!ht]
\centering
\subfigure[Density]{
\includegraphics[scale=.35]{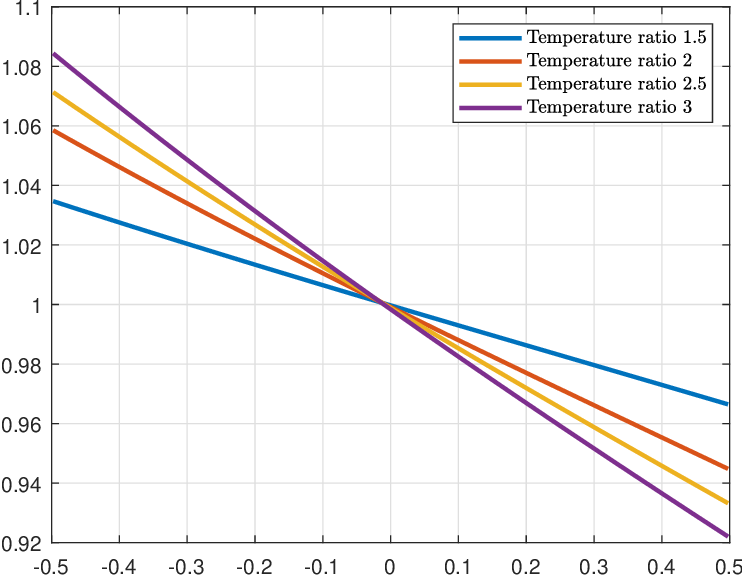}
} \quad
\subfigure[Temperature]{
\includegraphics[scale=.35]{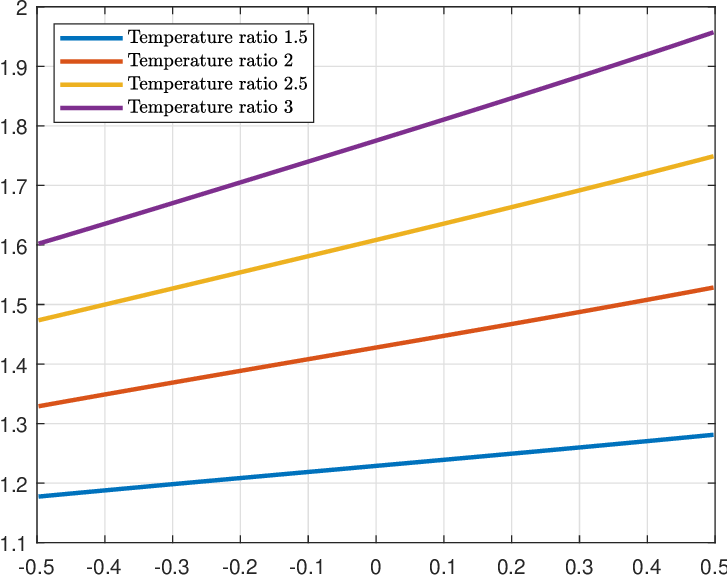}
}
\caption{Results of Fourier flows with $Kn = 1.0$.}
\label{fig:Fourier_tr}
\end{figure}

Our second example is the Fourier flow, where both walls are stationary but have different temperatures. In other words, we set $\bu_L = \bu_R = 0$ in \eqref{eq:maxwell_bc}, and $\theta_L$ and $\theta_R$ are different. In our experiments, we fix $\theta_L = 1$ and let $\theta_R$ take four values $1.5$, $2$, $2.5$, $3$. Since large temperatures indicate wider distribution functions, we choose $R = 3.5$ for $\theta_R = 1.5$ and $2$, and $R = 4$ for $\theta_R = 2.5$ and $3$. Some solutions are plotted in Figure \ref{fig:Fourier_Kn} and Figure \ref{fig:Fourier_tr}. Like the results of Couette flows, these figures show that the general trend of the solutions agrees with the expectation: both density and temperature curves are flatter for smaller temperature ratios or larger Knudsen numbers, and the temperature jump is stronger for greater $\Kn$.

\begin{figure}[!ht]
\centering
\subfigure[$Kn = 0.1$]{
\includegraphics[scale=.35]{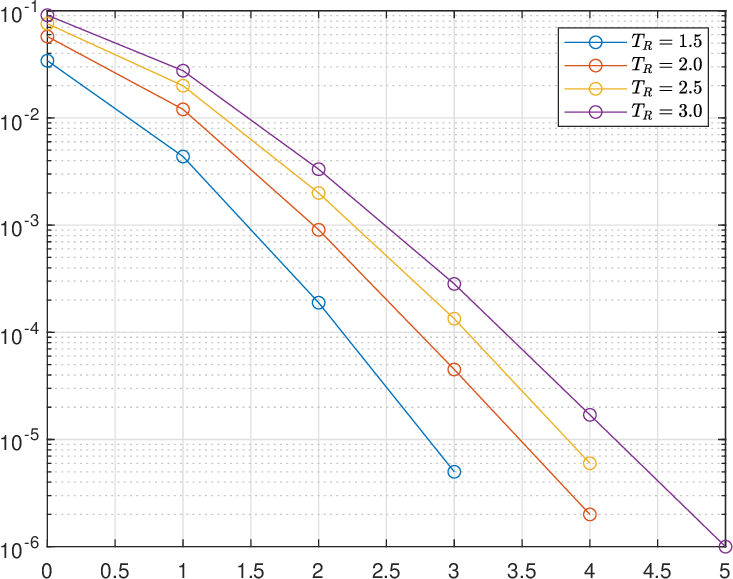}
} \quad
\subfigure[$Kn = 1$]{
\includegraphics[scale=.35]{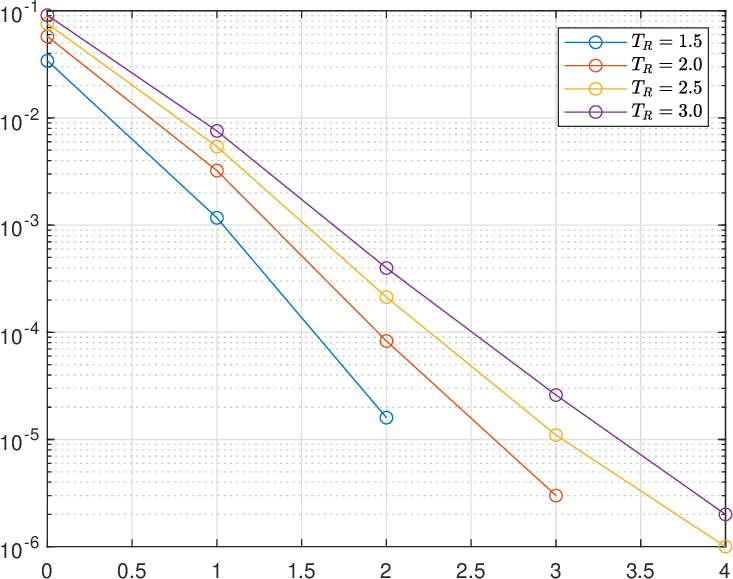}
}
\caption{Convergence of Newton iterations for Fourier flows. The horizontal axis denotes the number of iterations, and the vertical axis stands for the $L^2$ residual.}
\label{fig:Fourier_Newton}
\end{figure}

\begin{figure}[!ht]
\centering
\subfigure[$Kn = 0.1$]{
\includegraphics[scale=.35]{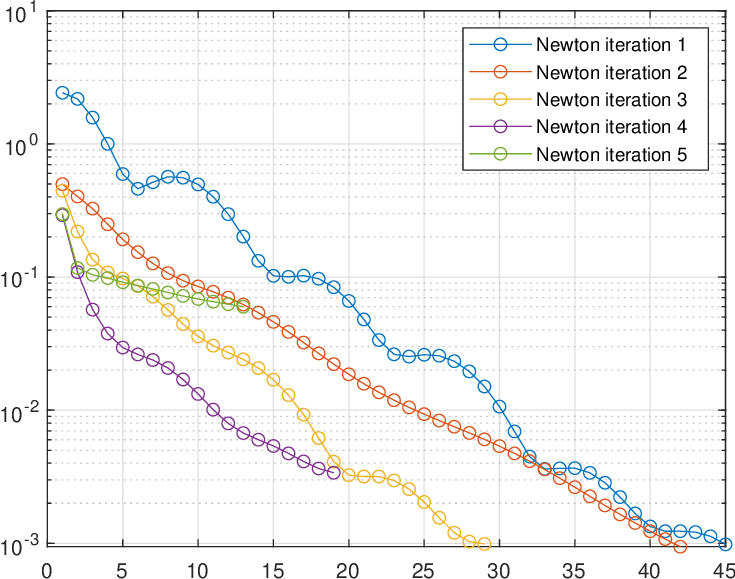}
} \quad
\subfigure[$Kn = 1$]{
\includegraphics[scale=.35]{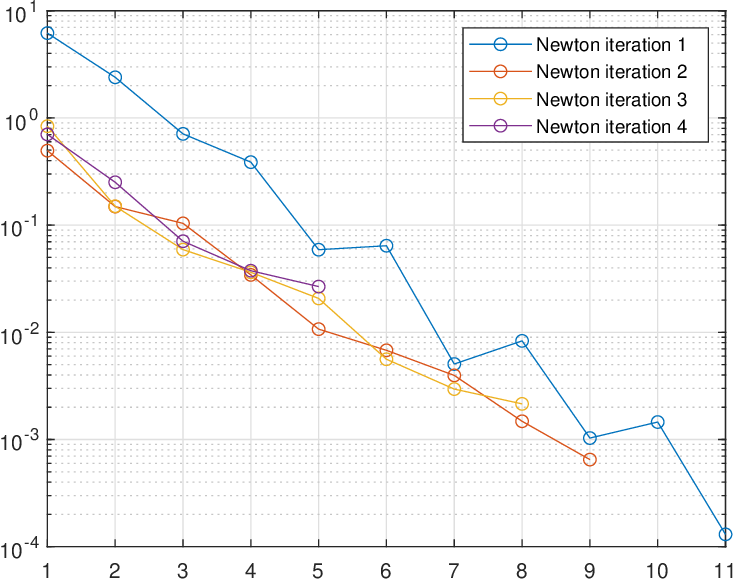}
}
\caption{Convergence of inner iterations for Fourier flows with $\theta_R = 3$. The horizontal axis denotes the number of iterations, and the vertical axis stands for the relative $L^2$ residual (the residual at the current step divided by the initial residual).}
\label{fig:Fourier_inner}
\end{figure}

Below, we again show some data to demonstrate the fast convergence of Newton's method. Table \ref{tab:Newton_Fourier} provides the number of Newton iterations for different parameters, and Figure \ref{fig:Fourier_Newton} gives more information on the error decay. The general behavior is the same as Couette flows: Newton's method takes only a few iterations to reach the desired residual, meaning that the binary collision operator is evaluated only a few times on each grid cell. To show some information on the inner iterations, here we pick $\theta_R = 3$ as an example and plot the error decay of inner iterations in Figure \ref{fig:Fourier_inner}. As expected, the case $\Kn = 0.1$ requires significantly larger numbers of iterations to achieve the relative residual $10^{-3}$. Note that here each iteration refers to one evaluation of the linearized collision operator on each grid cell.
\begin{table}[!ht]
\caption{Numbers of Newton iterations for Fourier flows}
\label{tab:Newton_Fourier}
\centering
\scalebox{.9}{%
\begin{tabular}{c|cccc}
  & $\theta_R = 1.5$ & $\theta_R = 2$ & $\theta_R = 2.5$ & $\theta_R = 3$ \\
  \hline
  $\Kn = 0.1$ & $3$ & $4$ & $4$ & $5$ \\
  $\Kn = 1$ & $3$ & $3$ & $4$ & $4$ \\
  $\Kn = 10$ & $2$ & $2$ & $3$ & $3$
\end{tabular}
}
\end{table}

In general, our tests verify that the modified Newton's method maintains the fast convergence rate, and effectively limits the number of evaluations of binary collision operators. Inner iterations, which involve linearized collision operators only, are to be improved in future works.

\section{Conclusion}
\label{sec:conclusion}

We have introduced an efficient numerical method with time complexity $O(JN^3\log N)$ (or $O(N^4 \log N)$ if $J = O(N)$) for the Boltzmann collision operator linearized about any local Maxwellian. However, straightforward application of the scheme may result in accuracy issues due to the catastrophic cancellations in Fourier transforms, and we introduced a suitable cutoff to the gain term to recover desired precision. Based on this algorithm, we modify the Newton method and achieve fast convergence rates when solving the steady-state Boltzmann equation.

The advantage of the modified Newton's method is that only one evaluation of the collision term per grid cell is needed in each iteration.
While some other methods such as the source iteration \cite{Fiveland1984}, the generalized synthetic iterative scheme \cite{Zhu2021}, and the symmetric Gauss-Seidel method \cite{Cai2025} may also have the same or a similar property, their convergence rates are generally slower than Newton's method, leading to a higher number of iterations, meaning more evaluations of the quadratic term.
Extra computational cost of our method comes from the solver of the linear system \eqref{eq:modified_lin_eq}, whose collision term requires a much lower computational cost using our method, and we are currently working on efficient numerical solvers of this equation.
Besides,
our ongoing work also includes extensions of this algorithm to more general gas molecules and generalization of the steady-state Boltzmann solver to multi-dimensional cases.

\bibliographystyle{siamplain}
\bibliography{Reference}
\end{document}


\maketitle

\section{Fourier transform of the gain term}
Here we present the detailed derivation for the Fourier transform of the gain term $\mL^+(\bv)$. Recall that
\begin{displaymath}
\mL^+(\bv) = \int_{\mathbb{R}^3} \int_{\mathbb{S}^2} \mB(|\bv - \bv_{\ast}|) [\mM(\bv') f(\bv_{\ast}') + f(\bv') \mM(\bv_{\ast}')] \,\dd \bsigma \,\dd \bv_{\ast}.
\end{displaymath}
The Fourier transform of the gain term can then be written as
\begin{displaymath}
\begin{aligned}
\widehat{\mL}^+(\bxi) &= 
 \int_{\mathbb{R}^3} \int_{\mathbb{R}^3} \int_{\mathbb{S}^2} \mathrm{e}^{-\image \bv \cdot \bxi} \mB(|\bv - \bv_{\ast}|) [\mM(\bv') f(\bv_{\ast}') + f(\bv') \mM(\bv_{\ast}')] \,\dd \bsigma \,\dd \bv_{\ast} \,\dd \bv \\
&= \int_{\mathbb{R}^3} \int_{\mathbb{R}^3} \int_{\mathbb{S}^2} \mathrm{e}^{-\image \bv \cdot \bxi} \mB(|\bv - \bv_{\ast}|) [r(\bv_{\ast}') + r(\bv')] \mM(\bv') \mM(\bv_{\ast}') \,\dd \bsigma \,\dd \bv_{\ast} \,\dd \bv \\
&= \frac{1}{(2\pi)^3}\int_{\mathbb{R}^3} \int_{\mathbb{R}^3} \int_{\mathbb{R}^3} \int_{\mathbb{S}^2}
\mathrm{e}^{-\image \bv \cdot \bxi} \mB(|\bv - \bv_{\ast}|) \hat{r}(\bdeta) \\
& \hspace{50pt} \times (\mathrm{e}^{\image \bv' \cdot \bdeta} + \mathrm{e}^{\image \bv_{\ast}' \cdot \bdeta})
\mM(\bv') \mM(\bv_{\ast}')  \,\dd \bsigma \,\dd \bv_{\ast} \,\dd \bv \,\dd \bdeta.
\end{aligned}
\end{displaymath}
Further simplification again requires the change of variables \eqref{eq:change_of_var}, which leads to
\begin{equation} \label{eq:MM_equality}
\mM(\bv') \mM(\bv_{\ast}') = \mM(\bv) \mM(\bv_{\ast}) =\mM^g(\bg) \mM^h(\bh),
\end{equation}
where
\begin{displaymath}
\mM^g(\bg) = \frac{\rho}{(4\pi\theta)^{3/2}} \exp\left( -\frac{|\bg|^2}{4\theta} \right), \qquad \mM^h(\bh) = \frac{\rho}{(\pi\theta)^{3/2}} \exp\left( -\frac{|\bh-\bu|^2}{\theta} \right),
\end{displaymath}
and the equality $\mM(\bv') \mM(\bv_{\ast}') = \mM(\bv) \mM(\bv_{\ast})$ can be derived easily from the momentum and energy conservation of the collision dynamics.
Hence,
\begin{displaymath}
\begin{aligned}
& \widehat{\mL}^+(\bxi) = \frac{1}{(2\pi)^3} \int_{\mathbb{R}^3} \int_{\mathbb{R}^3} \int_{\mathbb{R}^3} \int_{\mathbb{S}^2} \mB(|\bg|) \hat{r}(\bdeta) (\mathrm{e}^{\image |\bg|\bsigma \cdot \bdeta/2} + \mathrm{e}^{-\image |\bg| \bsigma \cdot \bdeta/2}) \\
& \hspace{50pt} \times \mathrm{e}^{-\image \bg \cdot \bxi/2} \mathrm{e}^{-\image \bh \cdot (\bxi-\bdeta)} \mM^g(\bg) \mM^h(\bh) \,\dd \bsigma \,\dd \bg \,\dd \bh \,\dd \bdeta.
\end{aligned}
\end{displaymath}
We can now integrate out $\bh$ and $\bsigma$ explicitly to get
\begin{displaymath}
\widehat{\mL}^+(\bxi) = \frac{1}{\pi^2} \int_{\mathbb{R}^3} \int_{\mathbb{R}^3} \mB(|\bg|) \hat{r}(\bdeta) \operatorname{sinc} \left( \frac{|\bg| |\bdeta|}{2} \right) \mathrm{e}^{-\image \bg \cdot \bxi/2} \mM^g(\bg) \widehat{\mM}^h(\bxi-\bdeta) \,\dd \bg \,\dd \bdeta,
\end{displaymath}
where $\widehat{\mM}^h(\cdot)$ denotes the Fourier transform of $\mM^h(\cdot)$ and can be calculated explicitly as
\begin{displaymath}
\widehat{\mM}^h(\bdeta) = \rho \exp (-\image \bu \cdot \bdeta - \theta |\bdeta|^2).
\end{displaymath}
Write $\bg$ in spherical coordinates $\bg = g \bomega$, allowing us to integrate out the angular variable $\bomega$:
\begin{equation} 
\begin{aligned}
\widehat{\mL}^+(\bxi) &= \frac{1}{\pi^2} \int_{\mathbb{R}^3} \int_0^{+\infty} \int_{\mathbb{S}^2} g^2 \mB(g) \hat{r}(\bdeta) \operatorname{sinc} \left( \frac{g |\bdeta|}{2} \right) \\
& \hspace{50pt} \times \mathrm{e}^{-\image g \bomega \cdot \bxi/2} \left[ \frac{\rho}{(4\pi\theta)^{3/2}} \exp \left( -\frac{g^2}{4\theta} \right) \right] \widehat{\mM}^h(\bxi-\bdeta) \,\dd \bomega \,\dd g \,\dd \bdeta \\
&= \frac{4}{\pi} \int_0^{+\infty} g^2 \mB(g) \left[ \frac{\rho}{(4\pi\theta)^{3/2}} \exp \left( -\frac{g^2}{4\theta} \right) \right] \operatorname{sinc}\left( \frac{g|\bxi|}{2}\right) \\
& \hspace{50pt} \times \left[\int_{\mathbb{R}^3} \operatorname{sinc}\left( \frac{g|\bdeta|}{2}\right) \hat{r}(\bdeta)\widehat{\mM}^h(\bxi-\bdeta) \,\dd\bdeta \right] \,\dd g.
\end{aligned}
\end{equation}
Note that this explicit integration with respect to $\bomega$ is made possible by introducing the new function $r(\bv)$, leading to a result that does not contain any integrals over $\mathbb{S}^2$.

However, this expression may fail if $f(\bv)$ does not decay as fast as $\mathcal{\mM}(\bv)$. In this case, the Fourier transform of $r(\bv)$ may not exist. This can be fixed by introducing the following function:
\begin{equation} \label{eq:hatH}
\hat{H}(g, \bxi) := \int_{\mathbb{R}^3} g^2 \operatorname{sinc}\left( \frac{g|\bdeta|}{2}\right) \hat{r}(\bdeta)\widehat{\mM}^h(\bxi-\bdeta) \,\dd\bdeta,
\end{equation}
so that
\begin{equation}
\widehat{\mL}^+(\bxi) = \frac{4}{\pi} \int_0^{+\infty} \mB(g) \left[ \frac{\rho}{(4\pi\theta)^{3/2}} \exp \left( -\frac{g^2}{4\theta} \right) \right] \operatorname{sinc}\left( \frac{g|\bxi|}{2}\right)  \hat{H}(g,\bxi) \,\dd g.
\end{equation}
We can apply the inverse Fourier transform to  \eqref{eq:hatH} to convert it to the physical space:
\begin{equation}
H(g, \bv) = \frac{1}{(2\pi)^3} \int_{\mathbb{R}^3} \hat{H}(g,\bxi) \mathrm{e}^{\image \bv \cdot \bxi} \,\mathrm{d}\bxi = 8\pi^2 \mM^h(\bv) \int_{\mathbb{R}^3} \frac{f(\bv-\bw)}{\mM(\bv-\bw)} \delta\left( |\bw| - \frac{g}{2} \right) \dd \bw,
\end{equation}
where the Dirac delta function comes from the inverse Fourier transform of the sinc function:
\begin{displaymath}
\frac{1}{(2\pi)^3} \int_{\mathbb{R}^3} \operatorname{sinc} \left( \frac{g|\bxi|}{2} \right) \mathrm{e}^{\image \bxi \cdot \bv} \,\mathrm{d}\bxi = \frac{1}{\pi g^2} \delta\left(|\bv| - \frac{g}{2} \right).
\end{displaymath}
Despite the Maxwellian in the denominator, the function $H(g,\bv)$ actually decays exponentially with respect $\bv$ once $f(\cdot)$ is bounded. In fact, if $|f(\bv)| \leqslant C$ for a certain positive constant $C$, straightforward calculation yields
\begin{displaymath}
\begin{aligned}
|H(g,\bv)| & \leqslant 8\pi^2 C \mM^h(\bv) \int_{\mathbb{R}^3} \frac{1}{\mM(\bv-\bw)} \delta \left(|\bw| - \frac{g}{2} \right) \dd \bw \\
& = 16\sqrt{2} \pi^3 C \frac{g\theta}{|\bv-\bu|}\exp\left( \frac{g^2}{8\theta} \right)\exp\left( -\frac{|\bv-\bu|^2}{2\theta} \right) \sinh \left( \frac{g|\bv-\bu|}{2\theta} \right).
\end{aligned}
\end{displaymath}
Thus, for a general bounded distribution function $f$, $\widehat{\mL}^+(\bxi)$ can be expressed as
\begin{displaymath}
\widehat{\mL}^+(\bxi) = \frac{4}{\pi} \int_0^{+\infty} \mB(g) \left[ \frac{\rho}{(4\pi\theta)^{3/2}} \exp \left( -\frac{g^2}{4\theta} \right) \right] \operatorname{sinc}\left( \frac{g|\bxi|}{2}\right) \int_{\mathbb{R}^3} H(g,\bv) \mathrm{e}^{-\image \bv \cdot \bxi} \,\dd\bxi \,\dd g,
\end{displaymath}
which agrees with \eqref{eq:ft_linear} in the main article.

\section{Fully discrete scheme for the steady-state Boltzmann equation}
In this section, we write down the scheme for the steady-state Boltzmann equation with discrete spatial and velocity variables.
Recall that the steady-state Boltzmann equation reads
\begin{displaymath}
v_1 \partial_x f = \frac{1}{\Kn} \mQ[f,f], \qquad x \in (x_L, x_R), \quad \bv \in \mathbb{R}^3,
\end{displaymath}
and the boundary conditions are given in \eqref{eq:maxwell_bc}.
As a preliminary study, we will only consider the first-order finite volume method for spatial discretization. The domain $[x_L, x_R]$ is decomposed uniformly into $N_x$ cells. The grid size is denoted by $\Delta x$, and the coordinate of the center of the $j$th cell is assumed to be $x_j$. The velocity discretization will follow the spectral method introduced in previous sections. Below, we will use $\bv_{\bk}$ to denote the collocation point in the velocity space with multi-index $\bk$, and use $f_{j,\bk}$ to denote the numerical approximation of the average of $f(x,\bv_{\bk})$ in the $j$th grid cell. The steady-state Boltzmann equation can then be discretized using the upwind method as
\begin{equation} \label{eq:discrete_eq}
v_{1,\bk}^+ \frac{f_{j,\bk} - f_{j-1,\bk}}{\Delta x} + v_{1,\bk}^- \frac{f_{j+1,\bk} - f_{j,\bk}}{\Delta x} = \frac{1}{\Kn} \mQ_{\bk}[f_{j,\cdot}, f_{j,\cdot}],
\end{equation}
where the collision term $\mQ_{\bk}[f_{j,\cdot}, f_{j,\cdot}]$ is obtained by applying the spectral method introduced in Section \ref{sec:fsm_binary} to the average distribution function in the $j$th grid cell and retrieving the point value of the collision term at $\bv_{\bk}$. For $j = 1$ and $j = N_x$, the above equation needs distribution functions in ghost cells, which are given according to \eqref{eq:maxwell_bc}. For instance,
\begin{displaymath}
f_{0,\bk} = \frac{1}{(2\pi\theta_L)^{3/2}} \left[ -\sqrt{\frac{2\pi}{\theta_L}} \left( \frac{L}{N} \right)^3 \sum_{\bk=-N}^N v_{1,\bk}^- f_{1,\bk} \right] \exp \left( -\frac{|\bv_{\bk} - \bu_L|^2}{2\theta_L} \right),
\end{displaymath}
and $f_{N_x+1,\bk}$ can be defined similarly. The condition \eqref{eq:total_mass} becomes the following constraint on the numerical solution:
\begin{displaymath}
\Delta x \left( \frac{L}{N} \right)^3 \sum_{j=1}^{N_x} \sum_{\bk = -N}^N f_{j,\bk} = C.
\end{displaymath}

Since our scheme is established based on Newton's method, we are required to have an initial value that is sufficiently close to the exact solution. In practice, this approximate solution can be obtained by various methods such as solving the Navier-Stokes-Fourier equations or the Boltzmann equation with a BGK-type collision term. In this work, we simply try to provide a reasonable initial guess:
\begin{equation} \label{eq:initial}
f_{j,\bk}^{(0)} = \frac{C}{x_R - x_L} \frac{1}{(2\pi \theta_j)^{3/2}} \exp \left( -\frac{|\bv_{\bk}|^2}{2\theta_j} \right),
\end{equation}
where
\begin{displaymath}
\theta_j = \frac{x_R - x_j}{x_R - x_L} \theta_L + \frac{x_j - x_L}{x_R - x_L} \theta_R.
\end{displaymath}
In \eqref{eq:initial}, the superscript ``$(0)$'' is the index of iteration. Our modified Newton's method updates the solution according to the following formula:
\begin{displaymath}
f_{j,\bk}^{(k+1)} = f_{j,\bk}^{(k)} - g_{j,\bk}^{(k)}.
\end{displaymath}
The quantity $g_{j,\bk}^{(k)}$ satisfies the linear system
\begin{equation} \label{eq:modified_Newton}
v_{1,\bk}^+ \frac{g_{j,\bk}^{(k)} - g_{j-1,\bk}^{(k)}}{\Delta x} + v_{1,\bk}^- \frac{g_{j+1,\bk}^{(k)} - g_{j,\bk}^{(k)}}{\Delta x} = \frac{1}{\Kn} \mL_{j,\bk}^{(k)}[g_{j,\cdot}^{(k)}] + r_{j,\bk}^{(k)},
\end{equation}
where
\begin{displaymath}
\mL_{j,\bk}^{(k)}[g_{j,\cdot}^{(k)}] = \mQ_{\bk}[f_{j,\cdot}^{(k)}, \mM_{j,\cdot}^{(k)}] + \mQ_{\bk}[\mM_{j,\cdot}^{(k)}, f_{j,\cdot}^{(k)}].
\end{displaymath}
and $r_{j,\bk}^{(k)}$ is the residual defined by
\begin{displaymath}
r_{j,\bk}^{(k)} = v_{1,\bk}^+ \frac{f_{j,\bk}^{(k)} - f_{j-1,\bk}^{(k)}}{\Delta x} + v_{1,\bk}^- \frac{f_{j+1,\bk}^{(k)} - f_{j,\bk}^{(k)}}{\Delta x} - \frac{1}{\Kn} \mQ_{\bk}[f_{j,\cdot}^{(k)}, f_{j,\cdot}^{(k)}].
\end{displaymath}
The boundary conditions of $g_{j,\bk}^{(k)}$ is the same as the boundary conditions of $f_{j,\bk}$.
Note that \eqref{eq:modified_Newton} is a discretization of \eqref{eq:lin_eq} using the upwind finite volume method.
Since $\mL_{j,\bk}^{(k)}[g_{j,\cdot}^{(k)}]$ can be evaluated efficiently, in each iteration,  the binary collision operator needs to be computed only once on each grid cell when calculating $r_{j,\cdot}^{(k)}$.

Various methods in the literature can be applied to solve the linear equation \eqref{eq:modified_Newton}, including the symmetric Gauss-Seidel method \cite{Cai2025} and the general synthetic iterative scheme \cite{SU_GSIS2020, Su2020}. These methods are highly efficient especially in the case where $\Kn$ is small. While these methods will be explored in our future works, here we will only apply a classical method known as source iteration to solve \eqref{eq:modified_Newton}, as mentioned in the main article.

The initial guess of $g_{j,\bk}^{(k)}$, denoted by $g_{j,\bk}^{(k,0)}$, is set to zero for all $j$ and $\bk$. Note that here the first superscript $k$ denotes the index of the Newton iteration, and the second superscript stands for the number of inner iterations to solve \eqref{eq:modified_Newton}. The source iteration updates the solution according to the following equation:
\begin{equation} \label{eq:source_iteration1}
\begin{aligned}
v_{1,\bk}^+ \frac{g_{j,\bk}^{(k,\ell+1)} - g_{j-1,\bk}^{(k,\ell+1)}}{\Delta x} + v_{1,\bk}^- \frac{g_{j+1,\bk}^{(k,\ell+1)} - g_{j,\bk}^{(k,\ell+1)}}{\Delta x} + \frac{\nu_j}{\Kn} g_{j,\bk}^{(k,\ell+1)} \qquad \\
= \frac{1}{\Kn} \left( \mL_{j,\bk}[g_{j,\cdot}^{(k,\ell)}] + \nu_j g_{j,\bk}^{(k,\ell)} \right) + r_{j,\bk}^{(k)}.
\end{aligned}
\end{equation}
The parameter $\nu_j$ is chosen such that $\mL_{j,\bk}[g_{j,\cdot}^{\{\ell\}}] + \nu_j g_{j,\bk}^{\{\ell\}}$ is positive for every $\bk$, which guarantees that the solution is nonnegative when $r_{j,\bk}$ is nonnegative for every $j$ and $\bk$. In our experiments, we do not strictly follow this rule and simply choose $\nu_j$ to be the density of $f_{j,\cdot}$ for both collision kernels in \eqref{eq:col_ker}.

In the stopping criterion \eqref{eq:res_r}, the $L^2$ norm of the residual is evaluated by
\begin{displaymath}
\sqrt{\Delta x \left( \frac{L}{N} \right)^3 \sum_{j=1}^{N_x} \sum_{\bk=-N}^{N-1} \left|r_{j,\bk}^{(k)}\right|^2}.
\end{displaymath}
The computation of residual in the inner iteration \eqref{eq:source_iteration1} is performed in a similar way.

\bibliographystyle{siamplain}
\bibliography{Reference}